\documentclass[aos,MSNbibl,seceqn,citesort,dvips]{arximspdf}
\usepackage{graphicx}

\doi{10.1214/11-AOS935}
\volume{39}
\issue{6}
\pubyear{2011}
\firstpage{3062}
\lastpage{3091}

\makeatletter

\newcommand{\fracb}[1]{\bigl(#1\bigr)}

\newtheorem{theorem}{Theorem}[section]

\newproclaim{rem}[theorem]{Remark}

\newcommand{\dconv}{\stackrel{\mathrm{d}}{\rightarrow}}
\newcommand{\pconv}{\stackrel{\mathrm{p}}{\rightarrow}}
\newcommand{\dequi}{\stackrel{\mathrm{d}}{=}}

\makeatother

\begin{document}
\begin{frontmatter}

\title{Unit roots in moving averages beyond first order\thanksref{T1}}
\runtitle{Unit roots in moving averages}

\thankstext{T1}{Supported in part by NSF Grants DMS-07-43459 and
DMS-11-07031.}

\begin{aug}
\author[A]{\fnms{Richard A.} \snm{Davis}\ead[label=e1]{rdavis@stat.columbia.edu}}
\and
\author[B]{\fnms{Li} \snm{Song}\corref{}\ead[label=e2]{li.song@barcap.com}}
\runauthor{R. A. Davis and L. Song}
\affiliation{Columbia University and Barclays Capital}
\address[A]{Department of Statistics\\
1255 Amsterdam Ave\\
Columbia University\\
New York, New York 10027\\
USA} %adresu isvedimo komanda gale!
\address[B]{Barclays Capital\\
745 7th Ave\\
New York, New York 10019\\
USA}
\end{aug}

% HISTORY:
\received{\smonth{6} \syear{2010}}
\revised{\smonth{7} \syear{2011}}

% ABSTRACT
%
\begin{abstract}
The asymptotic theory of various estimators based on Gaussian
likelihood has been developed for the unit root and near unit root
cases of a first-order moving average model. Previous studies of the
MA(1) unit root problem rely on the special autocovariance structure of
the MA(1) process, in which case, the eigenvalues and eigenvectors of
the covariance matrix of the data vector have known analytical forms.
In this paper, we take a different approach to first consider the joint
likelihood by including an augmented initial value as a parameter and
then recover the exact likelihood by integrating out the initial value.
This approach by-passes the difficulty of computing an explicit
decomposition of the covariance matrix and can be used to study unit
root behavior in moving averages beyond first order. The asymptotics of
the generalized likelihood ratio (GLR) statistic for testing unit roots
are also studied. The GLR test has operating characteristics that are
competitive with the locally best invariant unbiased (LBIU) test of
Tanaka for some local alternatives and dominates for all other
alternatives.
\end{abstract}

% KEYWORDS
%
\begin{keyword}[class=AMS]
\kwd{62M10}.
\end{keyword}
\begin{keyword}
\kwd{Unit roots}
\kwd{moving average}.
\end{keyword}

\end{frontmatter}

%s1 #&#
\section{Introduction}

In this paper we consider inference for moving average models that
possess one or more unit roots in the moving average polynomial. To
introduce the problem, let's first consider the MA(1) model given by
%
%e1.1 #&#
\begin{equation}\label{fmla1}
X_t=Z_t-\theta_0Z_{t-1},
\end{equation}
where $\theta_0\in\mathbb{R}$, $\{Z_t\}$ is a sequence of
independent and identically distributed (i.i.d.) random
variables with $\mathbf{E}Z_t=0, \mathbf{E}Z_t^2=\sigma_0^2$ and
density function $f_Z$. The MA(1) model is invertible if and only if
$|\theta_0|<1$,
since in this case $Z_t$ can be represented explicitly in terms of
past values of $X_t$, that is,
\[
Z_t=\sum_{j=0}^\infty\theta_0^jX_{t-j}.
\]
Under this invertibility constraint, standard estimation procedures
that produce asymptotically normal estimates are readily available.
For example, if~$\hat{\theta}$ represents the maximum likelihood
estimator, found by maximizing the Gaussian likelihood based on the
data $X_1,\ldots,X_n$, then it is well known (see Brockwell and
Davis \cite{bd}) that
%
%e1.2 #&#
\begin{equation}\label{f1}
\sqrt{n}(\hat{\theta}-\theta_0)\dconv
\mathrm{N}(0,1-\theta_0^2).
\end{equation}
From the form of the limiting variance in (\ref{f1}), the asymptotic
behavior of $\hat{\theta}$, let alone the scaling, is not
immediately clear in the unit root case corresponding to
$\theta_0=1$.

In the case $f_Z$ is Gaussian, the parameters $\theta_0$ and
$\sigma^2$ are not identifiab\-le without the constraint
$|\theta_0|\leq1$. In particular, the profile Gaussian
log-likelihood, obtained by concentrating out the variance
parameter, satisfies
%
%e1.3 #&#
\begin{equation}\label{sim}
L_n(\theta)=L_n(1/\theta).
\end{equation}
It follows that $\theta=1$ is a critical value of the profile
likelihood, and hence there is a positive probability that $\theta=1$
is indeed the maximum likelihood estimator. If $\theta_0=1$, then it
turns out that this probability does not vanish asymptotically (see,
e.g., Anderson and Takemura \cite{ander}, Tanaka
\cite{tanaka} and Davis and Dunsmuir \cite{dd}). This phenomenon is
referred to as the pile-up effect. For the case that $\theta_0=1$ or
is near one in the sense that $\theta_0=1+\gamma/n$, it was shown in
Davis and Dunsmuir \cite{dd} that
\[
n(\hat{\theta}-\theta_0)\dconv\xi_\gamma,
\]
where $\xi_\gamma$ is a random variable with a discrete component at
0, corresponding to the asymptotic pile-up effect, and a
continuous component. Most of the early work on this problem was
based on explicit knowledge of the eigenvectors and eigenvalues of the
covariance matrix for observations from an MA(1) process; see
Anderson and Takemura \cite{ander}. Recently,
Breidt et al. \cite{breidtd} and Davis and Song \cite{lsd1} looked at
model (\ref{fmla1}) under the Laplace likelihood and the
Gaussian likelihood without resorting to knowledge of the precise
form of eigenvectors and eigenvalues of the covariance matrix. Instead
they introduced an auxiliary variable, which acts like an initial value
and can be integrated out to form the likelihood.

With a couple exceptions, most of previous work dealt exclusively with
the zero-mean case. Sargan
and Bhargava \cite{sb} and Shephard \cite{neils} showed that for the
nonzero mean case, the so-called pile-up effect is more severe than
the zero mean case. Chen, Davis and Song \cite{cdra} extended the results
from Davis and Dunsmuir \cite{dd} to regression models with errors
from a noninvertible MA(1) process. It is shown that, with a mean
term present in the model, the pile-up probability goes up to more than 0.95.

The MA unit root problem can arise in many modeling contexts,
especially if a time series exhibits trend and seasonality.
For
example, in personal communication,\vadjust{\goodbreak} Richard Smith has mentioned the
presence of a unit in modeling some environmental time series related
to climate change \cite{smith}. After detrending and fitting an ARMA
model to the time series, Smith noticed that the MA component appeared
to have a unit root. One explanation for this phenomenon is that
detrending often involves the application of a high-pass filter to the
time series. In particular, the filter diminishes or obliterates any
power in the time series at low frequencies (including the 0
frequency). Consequently, the detrended data will have a spectrum with
0 power at frequency 0, which can only be fitted with ARMA process that
has a unit root in the MA component. While we only consider unit roots
in higher order moving averages in this paper, we believe the
techniques developed here will be applicable in a more general
framework of an ARMA model. This will be the subject of future investigation.

In this paper, we will use the stochastic approaches described in
\cite{breidtd} and~\cite{lsd1} to first study the case when there is
a regression component in the time series and errors are generated from
noninvertible MA(1). A vital issue in extending these results to
higher order MA models is the scaling required for the \textit{auxiliary
variable}. The scaling used for the regression problem in the MA(1)
case provides insight into the way in which the auxiliary variable
should be scaled in the higher order case. Quite surprisingly, when
there is
only one unit root in the MA(2) process, that is,
%
%e1.4 #&#
\begin{equation}\label{fmla16}
X_t=Z_t+c_1Z_{t-1}+c_2Z_{t-2},
\end{equation}
where $-c_1-c_2=1$ and $\{Z_t\}\sim$ i.i.d.
$(0,\sigma^2)$, the asymptotic distribution of the maximum likelihood estimator
$(\hat c_1,\hat c_2)'$ is exactly the same as in
invertible MA(2) case; see \cite{bd}. That is,
%
%e1.5 #&#
\begin{equation}\label{bbfm}
\sqrt{n}\pmatrix{
\hat c_1- c_1 \cr
\hat c_2- c_2}
\dconv\mathrm{N}\biggl(0,\left[\matrix{
1-c_2^2 & c_1(1-c_2) \cr
c_1(1-c_2) & 1-c_2^2}\right]\biggr).
\end{equation}
One difference, however, is that $\hat c_1$ and $\hat c_2$ are now totally
dependent asymptotically [$c_1(1-c_2)=(1-c_2)^2$].

As seen from (\ref{sim}), the first derivative of the profile likelihood
function is always 0 when $\theta=1$. Therefore, the development of
typical score tests or Wald tests is intractable in this case. Davis,
Chen and\vspace*{1pt}
Dunsmuir \cite{cdd} used the asymptotic result from \cite{dd} to
develop a test of $H_0\dvtx\theta=1$ based on~$\hat\theta_{\mathrm{MLE}}$ and the generalized likelihood ratio.
Interestingly, we will see that the estimator of the unit root in the
MA(2) case
has the same limit distribution as the corresponding estimator in the
MA(1) case. Thus, we can extend the methods used in the MA(1) case to
test for
unit roots in the MA(2) case.

The paper is organized as follows. In Section \ref{sec2}, we demonstrate our
method of proof applied to the MA(1) model with regression. This
case plays a key role in the extension to higher order MAs. Section
\ref{sec3} contains the results for the unit root problem in the MA(2) case.
In Section \ref{sec4}, we compare likelihood based\vadjust{\goodbreak} tests with Tanaka's
locally best invariant and unbiased (LBIU) test \cite{tpaper} for
testing the presence of a unit root. It is shown that the likelihood
ratio test performs quite well in comparison to the LBIU test. In
Section \ref{sec5}, numerical simulation results are presented to illustrate
the theory of Section \ref{sec3}. In Section \ref{sec6}, there is a brief discussion
that connects the auxiliary variables in higher order MAs with
terms in a regression model with MA(1) errors. Finally, in Section
\ref{sec7}, the procedure for handling the MA($q$) case with $q\geq3$ is
outlined. It is shown that the tools used in the MA(1) and MA(2)
cases are still applicable and are, in fact, sufficient in dealing
with higher order cases.

%s2 #&#
\section{MA(1) with nonzero mean}\label{sec2}

In this section, we will extend the methods of Breidt et al. \cite
{breidtd} and Davis and Song \cite{lsd1} to a regression model with
MA(1) errors. These results turn out to have connections with the
asymptotics in the higher order unit root cases (see Section \ref{sec6}). First,
consider the model
%
%e2.1 #&#
\begin{equation}\label{fmla13}
X_t=\sum_{k=0}^pb_{k0}f_k(t/n)+Z_t-\theta
_0Z_{t-1},
\end{equation}
where $\{Z_t\}$ is defined as in (\ref{fmla1}), $\theta_0=1$,
$b_{k0},k=0,\ldots,p$, are regression coefficients and
$f_k(t/n),k=0,\ldots,p$, are covariates at time $t$.
Notice that the covariates $f_{k}(t/n)$ are also assumed to be
functions on $[0,1]$. Note that the detrended series
$Y_t=X_t-\sum_{k=0}^pb_{k}f_k(t/n)$ has exactly the same
likelihood as the one for the zero-mean case. As shown in
\cite{lsd1}, by concentrating out the scale parameter $\sigma$,
maximizing the joint Gaussian likelihood is equivalent to minimizing
the following objective function:
%
%e2.2 #&#
\begin{equation}\label{fmla6}
l_n(\vec{b}, \theta,z_{\mathrm{init}}) = \sum_{t=0}^nz_t^2\qquad
\mbox{for
}|\theta|\leq1,
\end{equation}
where $\vec{b}=(b_0,\ldots,b_p)'$, $Z_{\mathrm{init}}=Z_0$, and $z_i$ is given by
\begin{eqnarray*}
z_i&=&Y_i+\theta Y_{i-1}+\cdots+\theta^{i-1}Y_1+\theta^iz_{\mathrm{init}}
\\
&=&\Biggl(X_i-\sum_{k=0}^pb_{k}f_k(i/n)\Biggr)+\theta
\Biggl(X_{i-1}-\sum_{k=0}^pb_{k}f_k\bigl((i-1)/n\bigr)
\Biggr)+\cdots\\
&&{}+\theta^{i-1}\Biggl(X_{1}-\sum_{k=0}^pb_{k}f_k
(1/n)\Biggr)+\theta^iz_{\mathrm{init}}\\
&=&\Biggl(Z_i-Z_{i-1}+\sum_{k=0}^pb_{k0}f_k(i/n)-\sum
_{k=0}^pb_{k}f_k(i/n)\Biggr)+\cdots\\
&&{}+\theta^{i-1}\Biggl(Z_{1}-Z_{0}+\sum_{k=0}^pb_{k0}f_k
(1/n)-\sum_{k=0}^pb_{k}f_k(1/n)\Biggr)+\theta
^iz_{\mathrm{init}}\\
&=&Z_i-(1-\theta)\sum_{j=0}^{i-1}\theta^{i-1-j}-\theta
^i(Z_0-z_{\mathrm{init}})\\
&&{}+\sum_{k=0}^p(b_{k0}-b_{k})\Biggl(\sum_{j=1}^i\theta^{i-j}f_k
(j/n)\Biggr)\\
:\!&=&Z_i-y_i+\sum_{k=0}^p(b_{k0}-b_{k})\Biggl(\sum_{j=1}^i\theta
^{i-j}f_k(j/n)\Biggr)\\
:\!&=&Z_i-w_i.
\end{eqnarray*}
As in \cite{lsd1}, we adopt the parametrization for $\theta$ and
$z_{\mathrm{init}}$ given by
\[
\theta=1+\frac{\beta}{n} \quad\mbox{and}\quad
z_{\mathrm{init}}=Z_0+\frac{\alpha\sigma_0}{\sqrt{n}}.
\]
Further set
%
%e2.3 #&#
\begin{equation}\label{rep1}
b_k=b_{k0}+\frac{\eta_k\sigma_0}{n^{{3/2}}}.
\end{equation}
Note that (\ref{rep1}) essentially characterizes the convergence
rate of the estimated $b_k$ to its true value $b_{k0}$. At first
glance, this parameterization may look odd since it depends on the
known parameter values, which are unavailable. This form of
reparameterization is used only for deriving the asymptotic theory
of the maximum likelihood estimators and not for estimation
purposes. One notes that $\beta=n(\theta-1)$,
$\eta_k=n^{3/2}(b_k-b_{k0})$, so that the asymptotics of the MLE
$\hat{\theta}$ and $\hat{b}_{k}$ of the associated parameters are
found by the limiting behavior of $\hat\beta=n(\hat\theta-1)$,
$\hat\eta_k=n^{3/2}(\hat b_k-b_{k0})$. Hence, it is not necessary to
know the true values in this analysis. The scaling $n^{3/2}$ for the
regression coefficients is an artifact of the assumption that the
regressors take the form $f_k(t/n)$ that is imposed on the problem.
This also results in a clean expression for the limit.

Under the\vspace*{1pt}
$(\vec{\eta},\beta,\alpha)$ parameterization, it is easily seen
\cite{lsd1}, minimizing~$l_n(\vec{b},\allowbreak\theta, z_{\mathrm{init}})$ with respect to
$\vec{b},\theta,z_{\mathrm{init}}$ is equivalent to minimizing the function
%
%e2.4 #&#
\begin{equation}
U_n(\vec{\eta},\beta,\alpha)\equiv\frac{1}{\sigma_0^2}
[l_n(\vec{b},\theta,z_{\mathrm{init}})-l_n(\vec{b}_0,1,Z_0)]
\end{equation}
with respect to $\vec{\eta},\beta$ and $\alpha$. Then using the weak
convergence results in Davis and Song \cite{lsd1},
\begin{eqnarray*}
&&
U_n(\vec{\eta},\beta,\alpha)\\[-3pt]
&&\qquad=\frac{1}{\sigma_0^2}\sum
_{i=0}^nz_i^2-Z_i^2=-2\sum_{i=0}^n\frac{w_iZ_i}{\sigma_0^2}+\sum
_{i=0}^n\frac{w_i^2}{\sigma_0^2}\\[-3pt]
&&\qquad\dconv2\beta\int_0^1\int_0^se^{\beta(s-t)}\,dW(t)\,dW(s)+2\alpha
\int_0^1e^{\beta s}\,dW(s)\\[-3pt]
&&\qquad\quad{}-2\sum_{k=0}^p\eta_k\int_0^1\biggl(\int_0^se^{\beta
(s-t)}f_k(t)\,dt\biggr)\,dW(s) \\[-3pt]
&&\qquad\quad{}+\int_0^1\Biggl(\beta\int_0^se^{\beta(s-t)}\,dW(t)+\alpha e^{\beta
s}-\sum_{k=0}^p\eta_k\int_0^se^{\beta(s-t)}f_k(t)\,dt\Biggr)^2\,ds\\[-3pt]
&&\qquad:=U(\vec{\eta},\beta,\alpha),\vspace*{-2pt}
\end{eqnarray*}
where ``$\dconv$'' indicates weak convergence on
$C(\mathbb{R}^{p+1}\times(-\infty,0]\times\mathbb{R})$. Throughout
this paper, when referring to convergence of stochastic processes on~$C(\mathbb{R}^{k})$, the notation ``$\dconv$'' (``$\pconv$'') means
convergence in distribution (probability) on~$C(\mathbb{K})$ where
$\mathbb{K}$ is any compact set in $\mathbb{R}^{k}$.

As a special case of a polynomial, set $f_k(t)=t^k$. In this case, the
limiting process $U(\vec{\eta},\beta,\alpha)$ is
\begin{eqnarray*}
U(\vec{\eta},\beta,\alpha)&=&2\beta\int_0^1\int_0^se^{\beta(s-t)}\,dW(t)\,dW(s)\\[-3pt]
&&+2\alpha\int_0^1e^{\beta s}\,dW(s)-2\sum_{k=0}^p\eta_k\int_0^1\biggl(\int_0^se^{\beta(s-t)}t^kdt\biggr)\,dW(s) \\[-3pt]
&&+\int_0^1\Biggl(\beta\int_0^se^{\beta(s-t)}\,dW(t)+\alpha e^{\beta s}-\sum_{k=0}^p\eta_k\int_0^se^{\beta(s-t)}t^k\,dt\Biggr)^2\,ds.\vspace*{-2pt}
\end{eqnarray*}
From now on we consider the simple case of just a nonzero mean,
that is, $p=0$ and $f_0(t)=1$. The formula further simplifies to
%
%e2.5 #&#
\begin{eqnarray}\label{fmla27}
U(\eta_0,\beta,\alpha)&=&2\beta\int_0^1\int_0^se^{\beta(s-t)}\,dW(t)\,dW(s)\nonumber
\\[-9.5pt]\\[-9.5pt]
&&+2\alpha\int_0^1e^{\beta s}\,dW(s)-2\eta_0\int_0^1\frac{1-e^{\beta s}}{\beta}\,dW(s)\nonumber\\[-3pt]
&&+\int_0^1\biggl(\beta\int_0^se^{\beta(s-t)}\,dW(t)+\alpha e^{\beta s}-\eta_0\frac{1-e^{\beta s}}{\beta}\biggr)^2\,ds.\nonumber\vspace*{-2pt}
\end{eqnarray}

As shown in \cite{lsd1}, one can recover the exact
likelihood by integrating out the initial parameter effects. More
specifically,
\begin{eqnarray*}
f(\mathbf{x}_n,z_{\mathrm{init}})&=&\prod_{t=0}^nf(z_t)\\[-3pt]
&=&\biggl(\frac{1}{\sqrt
{2\pi\sigma^2}}\biggr)^{n+1}\exp\biggl\{-\frac{\sum
_{t=0}^nz_t^2}{2\sigma^2}\biggr\}\\[-3pt]
&=&\biggl(\frac{1}{\sqrt{2\pi\sigma^2}}\biggr)^{n+1}\exp\biggl\{
-\frac{l_n(b_0,\theta,z_{\mathrm{init}})-l_n(b_{00},1,Z_0)+\sum
_{t=0}^nZ_t^2}{2\sigma^2}\biggr\}\\[-3pt]
&=&\biggl(\frac{1}{\sqrt{2\pi\sigma^2}}\biggr)^{n+1}\exp\biggl\{
-\frac{\sum_{t=0}^nZ_t^2}{2\sigma^2}\biggr\}\exp\biggl\{-\frac
{U_n(\eta_0,\beta,\alpha)\sigma_0^2}{2\sigma^2}\biggr\},
\end{eqnarray*}
integrating out the augmented variable $z_{\mathrm{init}}$ yields
%
%e2.6 #&#
\begin{eqnarray}\label{exact}
f(\mathbf{x}_n)&=&\int_{-\infty}^{+\infty}f(\mathbf
{x}_n,z_{\mathrm{init}})\,dz_{\mathrm{init}}\nonumber\\
&=&\biggl(\frac{1}{\sqrt{2\pi\sigma^2}}\biggr)^{n+1}\exp\biggl\{
-\frac{\sum_{t=0}^nZ_t^2}{2\sigma^2}\biggr\}\\
&&{}\times\frac{\sigma_0}{\sqrt
n}\int_{-\infty}^{+\infty}\exp\biggl\{-\frac{U_n(\eta_0,\beta
,\alpha)\sigma_0^2}{2\sigma^2}\biggr\}\,d\alpha.\nonumber
\end{eqnarray}

A similar argument as in \cite{lsd1} then shows that by profiling
out the variance parameter $\sigma^2$ the exact profile
log-likelihood $L_n(\eta_0,\beta)$ has the following property:
%
%e2.7 #&#
\begin{eqnarray}\label{fmlacon}
&&L_n(\eta_0,\beta)-L_n(\eta_0,0)\nonumber\\
&&\qquad\dconv L^*(\eta_0,\beta)\nonumber\\
&&\qquad=\log\int_{-\infty}^{+\infty}\exp\biggl\{-\frac{U(\eta_0,\beta
,\alpha)}{2}\biggr\}\,d\alpha\\
&&\qquad\quad{} -\log\int_{-\infty}^{+\infty}\exp\biggl\{-\frac{U(\eta
_0,0,\alpha)}{2}\biggr\}\,d\alpha.\nonumber
\end{eqnarray}
The weak convergence results on $C(\mathbb{R}^2)$ in (\ref{fmlacon})
can be used to show convergence in distribution of a sequence of
local maximizers of the objective functions $L_n$ to the maximizer
of the limit process $L$ provided the latter is unique almost
surely. This is the content of Remark 1 (see also Lemma 2.2) of
Davis, Knight and Liu~\cite{dkliu}, which for ease of reference, we
state a version here.\looseness=-1
%
%re2.1 #&#
\begin{rem}\label{remnew1}
Suppose $\{L_n(\cdot)\}$ is a sequence of stochastic processes which
converge in distribution to $L(\cdot)$ on $C(\mathbb{R}^k)$. If $L$ has
a unique maximizer $\tilde\beta$ a.s., then there exists a sequence of
local maximizers $\{\hat\beta_n\}$ of $\{L_n\}$ that converge in
distribution to $\tilde\beta$. Note that this is consistent with many
of the statements made in the classical theory for maximum likelihood
(see, e.g., Theorem 7.1.1 of Lehmann~\cite{ell}) and for inference in
nonstandard time series models; see Theorems 8.2.1 and 8.6.1 in
Rosenblatt~\cite{rosenblattm}, Breidt et al.~\cite{bdtrindade}, Andrews
et al. \cite{andrews} and Andrews et al.~\cite{acdavis}. In some cases,
for example, if the $\{L_n\}$ have concave sample paths, this can be
strengthened to convergence of the global\vadjust{\goodbreak} maximizers of $L_n$. See also
Davis, Chen and Dunsmuir~\cite{cdd}, Davis and Dunsmuir~\cite{dd97},
Breidt et a.l \cite{bdtrindade} for examples of other cases when
$\{L_n\}$ are not concave.
\end{rem}

Returning to our example, under the case when $\theta_0=1$, that is,
$\beta=0$, the limit of the exact likelihood is $L(\eta_0,\beta=0)$.
This corresponds to the situation of inference about the mean term
when it is known that the driving noise is an MA(1) process with a
unit root. Since the Gaussian likelihood is a quadratic function of
regression coefficients, $L(\eta_0,\beta=0)$ is a quadratic function
in $\eta_0$. Applying Remark \ref{remnew1}, we obtain that the MLE
$\hat{\eta}_0$ converges in distribution to $\tilde{\eta}_0$, the
global maximizer of $L(\eta_0,\beta=0)$. In particular,
$\tilde{\eta}_0$ is the value that makes
$\frac{\partial}{\partial\eta_0}L(\eta_0,\beta=0)=0$. Since
\begin{eqnarray*}
&&\frac{\partial}{\partial\eta_0}L(\eta_0,\beta=0)\\
&&\qquad=\frac{\int
_{-\infty}^{+\infty}\exp\{-{U(\eta_0,\beta=0,\alpha
)}/{2}\}(-({1}/{2})({\partial
U(\eta_0,\beta=0,\alpha)}/{\partial\eta_0}))\,d\alpha}{\int
_{-\infty}^{+\infty}\exp\{-{U(\eta_0,\beta=0,\alpha
)}/{2}\}\,d\alpha},
\end{eqnarray*}
where
\[
U(\eta_0,\beta=0,\alpha) = 2\alpha
W(1)-2\eta_0\int_0^1s\,dW(s)+\int_0^1(\alpha-s\eta_0)^2\,ds
\]
and
\[
\frac{\partial}{\partial\eta_0}U(\eta_0,\beta=0,\alpha) =
2\eta_0\int_0^1s^2\,ds-2\int_0^1s\,dW(s)-2\alpha\int_0^1s\,ds.
\]
Solving $\frac{\partial}{\partial\eta_0}L(\eta_0,\beta=0)=0$, we
find that
%
%e2.8 #&#
\begin{equation}
\tilde{\eta}_0=12\int_0^1s\,dW(s)-6W(1)\sim\mathrm{N}(0,12)
\end{equation}
and hence
%
%e2.9 #&#
\begin{equation}\label{fmla28}
n^{{3}/{2}}(\hat{b}_{0,n}-b_0)=\sigma_0\hat\eta_{0,n}\dconv
\mathrm{N}(0,12\sigma_0^2).
\end{equation}
This counter-intuitive result was also obtained earlier by Chen et
al. \cite{cdra}. It says the MLE of the mean term in the process
would behave like a normal distribution asymptotically, but with
convergence rate $n^{{3}/{2}}$. Notice that, even if one does
not know the true value of $\theta$, the MLE of the mean term would
still behave very much like (\ref{fmla28}) due to the large pile-up
effect in this case. However, the MLE is not asymptotically normal,
if both $b_0$ and $\theta$ are estimated.

%s3 #&#
\section{MA(2) with unit roots}\label{sec3}

The above approach, which also works in the invertible case, does not
rely on detailed knowledge of the form of the
eigenvectors and eigenvalues of the covariance matrix. Hence it has the
potential to work in higher order models\vadjust{\goodbreak} where the
eigenvector and eigenvalue structure is not known explicitly. We will
concentrate on the MA(2) process in this section and further
illustrate our methods.

\begin{figure}

\includegraphics{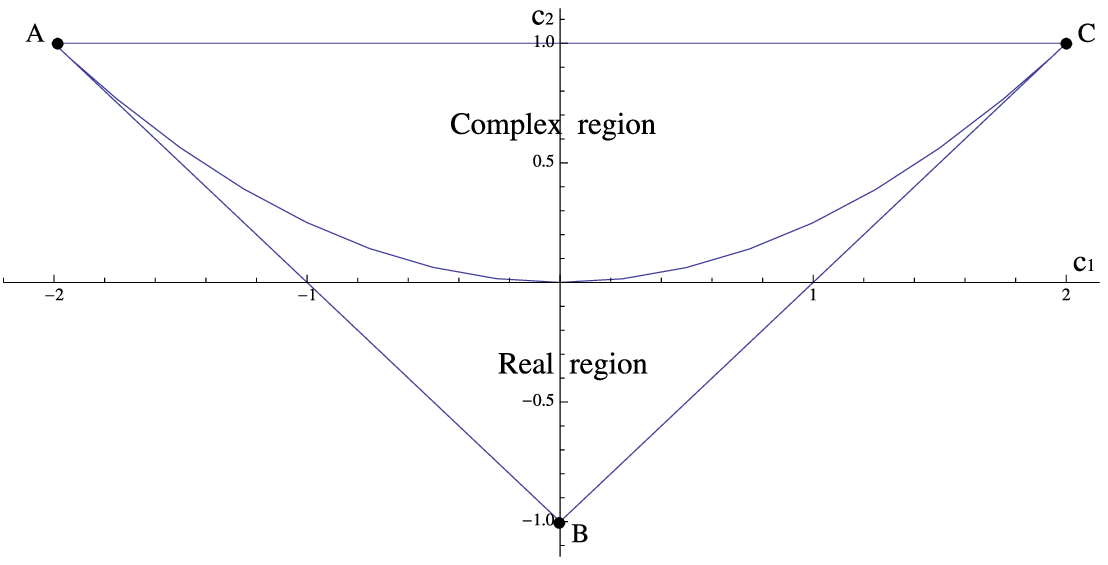}

\caption{$\bigtriangledown$ region defined by $-c_1-c_2\leq1,
c_1-c_2\leq1,
|c_2|\leq1$.}\label{figtri}
\end{figure}

In the following section, we consider the model given
in (\ref{fmla16}), where parameters
$c_1,c_2\in\bigtriangledown$, the triangular shaped region depicted
in Figure \ref{figtri}. The interior of this region corresponds to
the invertibility region of the parameter space. Note that the
triangular region is separated into complex roots and real roots of the
MA polynomial $1+c_1z+c_2z^2$ by a quadratic curve $c_1^2-4c_2=0$.

%f1 #&#

If the parameters are on the boundary of the $\bigtriangledown$
region, it indicates presence of unit roots. Otherwise, the model is
said to be invertible; see also Brockwell and Davis \cite{bd}. Model
(\ref{fmla16}) can also be represented in terms of the roots of the MA
polynomial by
\begin{eqnarray*}
X_t&=&(1+c_1\mathrm{B}+c_2\mathrm{B}^2)Z_t\\
&=&(1-\theta_0\mathrm{B})(1-\alpha_0\mathrm{B})Z_t,
\end{eqnarray*}
where $c_1=-\theta_0-\alpha_0$ and $c_2=\theta_0\alpha_0$.

%s3.1 #&#
\subsection{\texorpdfstring{Case 1: $|\alpha_0|<1$ and $\theta_0=1$}{Case 1: |alpha_0|<1 and theta_0=1}}

This case corresponds to the situation of only one unit root in the MA
polynomial, that is, the boundary AB in Figure \ref{figtri}. Let
$L_n(\theta,\alpha)$
be the profile likelihood of an MA(2) process. Again, we adopt the
parametrization
\[
\theta=1+\frac{\beta}{n}, \qquad\beta\leq0,
\]
and
\[
\alpha=\alpha
_0+\frac{\gamma}{\sqrt{n}}, \qquad\gamma\in\mathbb{R}.
\]
For convenience, define the intermediate process
$Y_t=(1-\alpha_0\mathrm{B})Z_t$ and observe that
\[
X_t=(1-\theta_0\mathrm{B})(1-\alpha_0\mathrm{B})Z_t=(1-\theta
_0\mathrm{B})Y_t.
\]
In the MA(2) case, two augmented initial variables
$Z_{\mathrm{init}}$ and $Y_{\mathrm{init}}$ are needed. These initial variables and
the joint likelihood have a simple form, that is,
%
%e3.1 #&#
\begin{equation}\label{fmla19}
Z_{\mathrm{init}}=Z_{-1} \quad\mbox{and}\quad
Y_{\mathrm{init}}=Z_0-\alpha_0Z_{\mathrm{init}},\vspace*{-3pt}
\end{equation}
\begin{eqnarray*}
f_{\mathbf{X},Y_{\mathrm{init}},Z_{\mathrm{init}}}(\mathbf
{x}_n,y_{\mathrm{init}},z_{\mathrm{init}})&=&f_{\mathbf{Y},Y_{\mathrm{init}},Z_{\mathrm{init}}}(\mathbf
{y}_n,y_{\mathrm{init}},z_{\mathrm{init}})\\
&=&f_{\mathbf{Z},Z_{\mathrm{init}}}(\mathbf{z}_n,z_{\mathrm{init}})\\
&=&\prod_{j=-1}^nf_Z(z_j).
\end{eqnarray*}
As what has been shown in the MA(1) case, the key of our method is
to calculate the formula for the residual $r_i:=Z_i-z_i$, which can
be obtained from
%
%e3.2 #&#
%e3.3 #&#
\begin{eqnarray}\quad
z_i&=&y_i+\alpha
y_{i-1}+\cdots+\alpha^{i-1}y_1+\alpha^iy_{\mathrm{init}}+\alpha
^{i+1}z_{\mathrm{init}}\nonumber\\
&=&\Biggl(\sum_{j=1}^i\theta^{i-j}X_j+\theta^iy_{\mathrm{init}}
\Biggr)+\alpha\Biggl(\sum_{j=1}^{i-1}\theta^{i-1-j}X_j+\theta
^{i-1}y_{\mathrm{init}}\Biggr)+\cdots\nonumber\\
&&{}+\alpha^{i-1}(X_1+\theta
y_{\mathrm{init}})+\alpha^iy_{\mathrm{init}}+\alpha^{i+1}z_{\mathrm{init}}\nonumber\\
&=&\sum_{j=1}^{i}\frac{\theta^{i-j+1}-\alpha^{i-j+1}}{\theta
-\alpha}X_j+\frac{\theta^{i+1}-\alpha^{i+1}}{\theta-\alpha
}y_{\mathrm{init}}+\alpha^{i+1}z_{\mathrm{init}}\nonumber\\
\label{fmla4th}
&=&Z_i-\frac{(\theta_0-\theta)(\theta-\alpha_0)}{\theta-\alpha
}\sum_{j=-1}^{i-1}\theta^{i-j-1}Z_j\\
&&{}-\frac{(\alpha_0-\alpha)(\theta_0-\alpha)}{\theta-\alpha}\sum
_{j=-1}^{i-1}\alpha^{i-j-1}Z_j+\frac{\theta^{i+1}-\alpha
^{i+1}}{\theta-\alpha}(y_{\mathrm{init}}-Y_0)\nonumber\\
&&{}+\alpha^{i+1}(z_{\mathrm{init}}-Z_{-1})+(\theta_0-\theta)\frac{\theta
^{i+1}-\alpha^{i+1}}{\theta-\alpha}Z_{-1}\nonumber\\
\label{fmla20}
&=&Z_i-r_i,
\end{eqnarray}
where the fourth equation (\ref{fmla4th}) comes from the fact that
$X_j=Z_j-(\theta_0+\alpha_0)\times Z_{j-1}+\theta_0\alpha_0Z_{j-2}$ and
$Y_0=Z_0-\alpha_0Z_{-1}$. Therefore, the residuals $r_i$ are given
by
%
%e3.4 #&#
\begin{eqnarray}\label{resi}
r_i&=&\frac{(\theta_0-\theta)(\theta-\alpha_0)}{\theta-\alpha
}\sum_{j=-1}^{i-1}\theta^{i-j-1}Z_j\nonumber\\
&&{}+\frac{(\alpha_0-\alpha)(\theta_0-\alpha)}{\theta-\alpha}\sum
_{j=-1}^{i-1}\alpha^{i-j-1}Z_j-\frac{\theta^{i+1}-\alpha
^{i+1}}{\theta-\alpha}(y_{\mathrm{init}}-Y_0)\\
&&{}-\alpha^{i+1}(z_{\mathrm{init}}-Z_{-1})-(\theta_0-\theta)\frac{\theta
^{i+1}-\alpha^{i+1}}{\theta-\alpha}Z_{-1}.\nonumber
\end{eqnarray}
Notice that the residuals $r_i$ no longer have a neat form as in the
MA(1) case. This is what makes the MA(2) case more
interesting yet more complicated.

In the following calculations, let
\[
y_{\mathrm{init}}=Y_0+\frac{\sigma_0\eta_1}{\sqrt{n}}
\quad\mbox{and}\quad
z_{\mathrm{init}}=Z_{-1}+\frac{\sigma_0\eta_2}{\sqrt{n}}.
\]
With a similar argument as in \cite{lsd1}, we opt to minimize the
objective function
%
%e3.5 #&#
\begin{equation}\label{obj}
U_n(\beta,\gamma,\eta_1,\eta_2)=-2\sum_{i=-1}^n\frac
{r_iZ_i}{\sigma_0^2}+\sum_{i=-1}^n\frac{r_i^2}{\sigma_0^2}.
\end{equation}
First note that $r_i=A_i+B_i+C_i+D_i$, where
\begin{eqnarray*}
A_i&:=&\frac{(\theta_0-\theta)(\theta-\alpha_0)}{\theta-\alpha
}\sum_{j=-1}^{i-1}\theta^{i-j-1}Z_j-\frac{\theta^{i+1}-\alpha
^{i+1}}{\theta-\alpha}(y_{\mathrm{init}}-Y_0),\\
B_i&:=&\frac{(\alpha_0-\alpha)(\theta_0-\alpha)}{\theta-\alpha
}\sum_{j=-1}^{i-1}\alpha^{i-j-1}Z_j,\\
C_i&:=&-\alpha^{i+1}(z_{\mathrm{init}}-Z_{-1}),\\
D_i&:=&-(\theta_0-\theta)\frac{\theta^{i+1}-\alpha^{i+1}}{\theta
-\alpha}Z_{-1}.
\end{eqnarray*}
To determine the weak limit of
$-2\sum_{i=-1}^n\frac{r_iZ_i}{\sigma_0^2}$ in (\ref{obj}) in the
continuous function space, note that
%
%e3.6 #&#
\begin{eqnarray}\label{firstpart}
-2\sum_{i=-1}^n\frac{A_iZ_i}{\sigma_0^2}&=&2\frac{(\theta-\theta
_0)(\theta-\alpha_0)}{\theta-\alpha}\sum_{i=-1}^n\sum
_{j=-1}^{i-1}\theta^{i-j-1}\frac{Z_j}{\sigma_0}\frac{Z_i}{\sigma
_0}\nonumber\\
&&{}+\frac{2\eta_1}{\sqrt{n}(\theta-\alpha)}\sum_{i=-1}^n\theta
^{i+1}\frac{Z_i}{\sigma_0}-\frac{2\eta_1}{\sqrt{n}(\theta-\alpha
)}\sum_{i=-1}^n\alpha^{i+1}\frac{Z_i}{\sigma_0}\nonumber\\
&=&2\frac{\beta(1-\alpha_0+{\beta}/{n})}{1-\alpha_0+
{\beta}/{n}-{\gamma}/{\sqrt{n}}}\sum_{i=-1}^n\sum
_{j=-1}^{i-1}\biggl(1+\frac{\beta}{n}\biggr)^{i-j-1}\frac
{Z_j}{\sigma_0}\frac{Z_i}{\sqrt{n}\sigma_0}\nonumber\\[-8pt]\\[-8pt]
&&{}+\frac{2\eta_1}{(1-\alpha_0+{\beta}/{n}-{\gamma
}/{\sqrt{n}})}\sum_{i=-1}^n\biggl(1+\frac{\beta}{n}
\biggr)^{i+1}\frac{Z_i}{\sqrt{n}\sigma_0}\nonumber\\
&&{}-\frac{2\eta_1}{\sqrt{n}(1-\alpha_0+{\beta}/{n}-
{\gamma}/{\sqrt{n}})}\sum_{i=-1}^n\biggl(\alpha_0+\frac{\gamma
}{\sqrt{n}}\biggr)^{i+1}\frac{Z_i}{\sqrt{n}\sigma_0}\nonumber\\
&\dconv&2\beta\int_0^1\int_0^se^{\beta(s-t)}\,dW(t)\,dW(s)+\frac
{2\eta_1}{1-\alpha_0}\int_0^1e^{\beta
s}\,dW(s),\nonumber
\end{eqnarray}
where the last term disappears in the limit due to the fact
that $|\alpha_0|<1$. Similarly, we have
%
%e3.7 #&#
%e3.8 #&#
\begin{eqnarray}
-2\sum_{i=-1}^n\frac{B_iZ_i}{\sigma_0^2}&=&2\frac{(\alpha-\alpha
_0)(1-\alpha)}{\theta-\alpha}\sum_{i=-1}^n\sum_{j=-1}^{i-1}\alpha
^{i-j-1}\frac{Z_j}{\sigma_0}\frac{Z_i}{\sigma_0}\nonumber\\
&=&2\frac{\gamma(1-\alpha_0-{\gamma}/{\sqrt{n}})}{1-\alpha
_0+{\beta}/{n}-{\gamma}/{\sqrt{n}}}\nonumber\\
&&{}\times\sum_{i=-1}^n\sum
_{j=-1}^{i-1}\biggl(\alpha_0+\frac{\gamma}{\sqrt{n}}
\biggr)^{i-j-1}\frac{Z_j}{\sigma_0}\frac{Z_i}{\sqrt{n}\sigma_0}\nonumber
\\
\label{fmla17}
&=&2\gamma\sum_{i=-1}^n\sum_{j=-1}^{i-1}\alpha_0^{i-j-1}\frac
{Z_j}{\sigma_0}\frac{Z_i}{\sqrt{n}\sigma_0}+o_p(1)\\
\label{fmla18}
&\dconv&2\gamma N,
\end{eqnarray}
where $N\sim\mathrm{N}(0,\frac{1}{1-\alpha_0^2})$. The third
equality holds because
$|\alpha_0|$ is strictly smaller than 1, and $o_p(1)$ is uniform in
$\gamma$ on any compact set of $\mathbb{R}$. The weak convergence
from (\ref{fmla17}) to (\ref{fmla18}) follows from martingale
central limit theorem; see Hall and Heyde \cite{heyde}. It can also be
shown that $N$ and the $W(t)$ process from (\ref{firstpart}) are
independent; see Theorem 2.2 in Chan and Wei \cite{weichan}.

Following similar arguments, it is easy to show that
\[
-2\sum_{i=-1}^n\frac{C_iZ_i}{\sigma_0^2}\pconv0 \quad\mbox
{and}\quad -2\sum_{i=-1}^n\frac{D_iZ_i}{\sigma_0^2}\pconv0.
\]
For the second term in (\ref{obj}), writing
\begin{eqnarray*}
\sum_{i=-1}^n\frac{r_i^2}{\sigma_0^2}&=&\sum_{i=-1}^n\frac
{A_i^2+B_i^2+C_i^2+D_i^2}{\sigma_0^2}\\
&&{}+\sum_{i=-1}^n\frac
{2A_iB_i+2A_iC_i+2A_iD_i+2B_iC_i+2B_iD_i+2C_iD_i}{\sigma_0^2},
\end{eqnarray*}
and using Corollary 2.10 in \cite{lsd1}, we have
%
%e3.9 #&#
%e3.10 #&#
\begin{eqnarray}
\sum_{i=-1}^n\frac{A_i^2}{\sigma_0^2}&\dconv&\int_0^1\biggl(\beta
\int_0^se^{\beta(s-t)}\,dW(t)+\frac{\eta_1}{1-\alpha_0}e^{\beta
s}\biggr)^2\,ds,\\
\sum_{i=-1}^n\frac{B_i^2}{\sigma_0^2}&\pconv&\gamma^2\operatorname{var}(N).
\end{eqnarray}
Moreover, it is relatively easy to show that
%
%e3.11 #&#
\begin{equation}
\sum_{i=-1}^n\frac{C_i^2}{\sigma_0^2}\pconv0 \quad\mbox{and}\quad
\sum_{i=-1}^n\frac{D_i^2}{\sigma_0^2}\pconv0.
\end{equation}
Next we show that all the cross product terms also vanish in the
limit, namely,
%
%e3.12 #&#
\begin{equation}\label{fmla24}\quad
\sum_{i=-1}^n\frac
{2A_iB_i+2A_iC_i+2A_iD_i+2B_iC_i+2B_iD_i+2C_iD_i}{\sigma_0^2}\pconv0.
\end{equation}
Here we only give the details for showing
$\sum_{i=-1}^n\frac{A_iB_i}{\sigma_0^2}\pconv0$; the other cases can
be proved in an analogous manner. Notice that for any fixed $M>0$
and any $\beta\in[-M,0]$,
%
%e3.13 #&#
\begin{eqnarray}\label{fmla23}
\sum_{i=-1}^n\frac{A_iB_i}{\sigma_0^2}&=&\frac{{\beta
}/{n}(1-\alpha_0+{\beta}/{n}){\gamma}/{\sqrt
{n}}(1-\alpha_0-{\gamma}/{\sqrt{n}})}{
(1-\alpha_0+{\beta}/{n}-{\gamma}/{\sqrt{n}}
)^2}\nonumber\\
&&{}\times\sum_{i=-1}^n\Biggl(\sum_{j=-1}^{i-1}\alpha^{i-j-1}\frac
{Z_j}{\sigma_0}\Biggr)\Biggl(\sum_{j=-1}^{i-1}\theta^{i-j-1}\frac
{Z_j}{\sigma_0}\Biggr)\nonumber\\
&&{}+\frac{({\gamma}/{\sqrt{n}})(1-\alpha_0-{\gamma
}/{\sqrt{n}})({\eta_1}/{\sqrt{n}})}{(1-\alpha
_0+{\beta}/{n}-{\gamma}/{\sqrt{n}})^2}\nonumber\\
&&\hspace*{11pt}{}\times\sum
_{i=-1}^n\Biggl[(\theta^{i+1}-\alpha^{i+1})\sum_{j=-1}^{i-1}\alpha
^{i-j-1}\frac{Z_j}{\sigma_0}\Biggr]\\
&=&\frac{\beta\gamma}{n}\sum_{i=-1}^n\Biggl(\sum
_{j=-1}^{i-1}\alpha_0^{i-j-1}\frac{Z_j}{\sigma_0}\Biggr)\Biggl(\sum
_{j=-1}^{i-1}\biggl(1+\frac{\beta}{n}\biggr)^{i-j-1}\frac
{Z_j}{\sqrt{n}\sigma_0}\Biggr)\nonumber\\
&&{}+\frac{\gamma\eta_1}{n}\sum_{i=1}^n\Biggl[\biggl(1+\frac{\beta
}{n}\biggr)^{i+1}\sum_{j=-1}^{i-1}\alpha_0^{i-j-1}\frac{Z_j}{\sigma
_0}\Biggr]\nonumber\\
&&\hspace*{11pt}{}\times\frac{\gamma\eta_1}{n}\sum_{i=-1}^n\sum_{j=-1}^{i-1}\alpha
_0^{2i-j}\frac{Z_j}{\sigma_0}+o_p(1),\nonumber
\end{eqnarray}
where $o_p(1)$ is uniform in $\beta$ and $\gamma$ on any compact set
in $\mathbb{R}^-\times\mathbb{R}$. Setting
$R_{i}=\sum_{j=-1}^i\alpha_0^{i-j}Z_j/\sigma_0$, it follows that
$R_i$ is a stationary AR(1) process satisfying
\[
R_i=\alpha_0R_{i-1}+Z_i/\sigma_0.
\]
Since $|\alpha_0|<1$, we can apply Theorem 3.7. in Tanaka
\cite{tanaka} to obtain
\[
S_n(t):=\frac{1}{\sqrt{n}}\sum_{i=0}^{[nt]}R_i\dconv\tilde{\alpha}S(t),\vspace*{-2pt}
\]
where\vspace*{1pt}
$\tilde{\alpha}=\sum_{l=0}^\infty\alpha_0^l=\frac{1}{1-\alpha_0}$
and $S(t)$ is a standard Brownian motion. Also, since $R_i$ is
adapted to the $\sigma$-fields $\mathcal{F}_i$ generated by
$Z_0,\ldots, Z_i$. By Theorem~2.1 in \cite{lsd1}, we obtain
\[
\frac{1}{\sqrt{n}}\sum_{i=-1}^n\biggl(1+\frac{\beta}{n}
\biggr)^{i+1}R_{i-1}\dconv\tilde{\alpha}\int_0^1e^{\beta
s}\,dS(s)\qquad \mbox{on }C[-M,0].\vspace*{-2pt}
\]
Therefore,
%
%e3.14 #&#
\begin{equation}\label{aux1}\quad
\frac{\gamma\eta_1}{n}\sum_{i=-1}^n\Biggl[\biggl(1+\frac{\beta
}{n}\biggr)^{i+1}\sum_{j=-1}^{i-1}\alpha_0^{i-j-1}\frac{Z_j}{\sigma
_0}\Biggr]\pconv0\qquad
\mbox{on }C[-M,0].\vspace*{-2pt}
\end{equation}
It is also easy to see that
%
%e3.15 #&#
\begin{equation}\label{aux6}
\frac{\gamma\eta_1}{n}\sum_{i=-1}^n\sum_{j=-1}^{i-1}\alpha
_0^{2i-j}\frac{Z_j}{\sigma_0}=\frac{\gamma\eta_1}{n}\sum
_{i=1}^n\alpha_0^{i+1}R_{i-1}\pconv0.\vspace*{-2pt}
\end{equation}
Since
%
%e3.16 #&#
\begin{equation}\label{aux5}
\sum_{i=-1}^n\Biggl(\sum_{j=-1}^{i-1}\biggl(1+\frac{\beta}{n}
\biggr)^{i-j-1}\frac{Z_j}{\sqrt{n}\sigma_0}\Biggr)\frac{R_{i-1}}{\sqrt
{n}}\vspace*{-2pt}
\end{equation}
is in the form of the double sum in Theorem 2.8 in \cite{lsd1},
except that $\{R_{i}\}$ is no longer a martingale difference
sequence. However, we can still follow the proof of Theorem 2.8 in
\cite{lsd1} and show that (\ref{aux5}) has a nondegenerate weak
limit in $C[-M,0]$. It follows that
%
%e3.17 #&#
\begin{eqnarray}\label{aux2}
&&\frac{\beta\gamma}{n}\sum_{i=-1}^n\Biggl(\sum_{j=-1}^{i-1}\alpha
_0^{i-j-1}\frac{Z_j}{\sigma_0}\Biggr)\Biggl(\sum_{j=-1}^{i-1}
\biggl(1+\frac{\beta}{n}\biggr)^{i-j-1}\frac{Z_j}{\sqrt{n}\sigma
_0}\Biggr)\nonumber\\[-9pt]\\[-9pt]
&&\qquad=\frac{\beta\gamma}{\sqrt{n}}\sum_{i=-1}^n\Biggl(\sum
_{j=-1}^{i-1}\biggl(1+\frac{\beta}{n}\biggr)^{i-j-1}\frac
{Z_j}{\sqrt{n}\sigma_0}\Biggr)\frac{R_{i-1}}{\sqrt{n}}\pconv
0.\nonumber\vspace*{-2pt}
\end{eqnarray}
Thus, combining (\ref{aux1}), (\ref{aux6}) and (\ref{aux2}), we
conclude that the terms in~(\ref{fmla23}) go to 0 in probability on
$C[-M,0]$. The convergence in probability of the other terms in (\ref
{fmla24}) can also be proved in a similar way. To sum up, we have shown
the key stochastic process convergence
result, that is,
%
%e3.18 #&#
\begin{eqnarray}\label{fmla25}\quad
&&U_n(\beta,\gamma,\eta_1,\eta_2)\nonumber\\[-2pt]
&&\qquad\dconv U(\beta,\gamma,\eta_1)\nonumber\\[-9pt]\\[-9pt]
&&\qquad=2\beta\int_0^1\int_0^se^{\beta(s-t)}\,dW(t)\,dW(s)+2\gamma
N+\frac{2\eta_1}{1-\alpha_0}\int_0^1e^{\beta s}\,dW(s)\nonumber\\[-2pt]
&&\qquad\quad{}+\int_0^1\biggl(\beta\int_0^se^{\beta(s-t)}\,dW(t)+\frac{\eta
_1}{1-\alpha_0}e^{\beta
s}\biggr)^2\,ds+\gamma^2\operatorname{var}(N).\nonumber\vspace*{-2pt}
\end{eqnarray}

Using (\ref{fmla25}), one can easily derive the
asymptotics for the exact profile log-likelihood denoted by
$L_n(\beta,\gamma)$. In particular,
%
%e3.19 #&#
%e3.20 #&#
\begin{eqnarray}
&&L_n(\beta,\gamma)-L_n(0,0)\nonumber\\
\label{fmla30}
&&\qquad\dconv\log\int_{-\infty}^{+\infty}\exp\biggl\{-\frac{U(\beta
,\gamma,\eta_1)}{2}\biggr\}\,d\eta_1\\
&&\qquad\quad{}-\log\int_{-\infty}^{+\infty}\exp\biggl\{-\frac{U(0,0,\eta
_1)}{2}\biggr\}\,d\eta_1\nonumber\\
&&\qquad:=L^*(\beta,\gamma)\nonumber\\
\label{fmla26}
&&\qquad=-\gamma
N-\frac{\gamma^2}{2}\operatorname{var}(N)+\log\int_{-\infty}^{+\infty
}\exp\biggl\{-\frac{U(\beta,\eta^*)}{2}\biggr\}\,d\eta^*\\
&&\qquad\quad{}-\log\int_{-\infty}^{+\infty}\exp\biggl\{-\frac{U(0,\eta
^*)}{2}\biggr\}\,d\eta^*,\nonumber
\end{eqnarray}
where $\eta^*=\frac{\eta_1}{1-\alpha_0}$ and $U(\beta,\eta^*)$ is
given by
%
%e3.21 #&#
\begin{eqnarray}\label{fmla5}
U(\beta,\eta^*)&=&2\int_0^1\biggl[\beta\int_0^se^{\beta
(s-t)}\,dW(t)+\eta^*
e^{\beta s}\biggr]\,dW(s)\nonumber\\[-8pt]\\[-8pt]
&&{}+\int_0^1\biggl[\beta\int_0^se^{\beta(s-t)}\,dW(t)+\eta^* e^{\beta
s}\biggr]^2\,ds,\nonumber
\end{eqnarray}
which is the limiting process of the joint likelihood obtained in
the unit root MA(1) case, see also Davis and Song \cite{lsd1}. We
state the key result of this paper in the following theorem.
%
%th3.1 #&#
\begin{theorem}\label{thm4}
Consider the model given in (\ref{fmla16}) with two roots $\theta$
and~$\alpha$ which are parameterized by
\[
\theta=1+\frac{\beta}{n}  \quad\mbox{and}\quad \alpha=\alpha
_0+\frac{\gamma}{\sqrt{n}}.
\]
Denote the profile log-likelihood based on a Gaussian likelihood as
$L_n(\beta,\gamma)$. Then $L_n(\beta,\gamma)$ satisfies
\[
L_n(\beta,\gamma)-L_n(0,0)\dconv L^*(\beta,\gamma) \qquad
\mbox{on
}C([-\infty,0]\times\mathbb{R}),
\]
where
%
%e3.22 #&#
\begin{eqnarray}\label{newadded}
L^*(\beta,\gamma)&=&-\gamma N-\frac{\gamma^2}{2}
\operatorname{var}(N)+U^*(\beta)\nonumber\\[-9pt]\\[-9pt]
&\dequi&-\gamma
N-\frac{\gamma^2}{2}\operatorname{var}(N)+\frac{1}{2}Z_0(\beta
).\nonumber
\end{eqnarray}
The processes
$U^*(\beta)$ and $Z_0(\beta)$ are defined by
%
%e3.23 #&#
\begin{eqnarray}\label{defu}
U^*(\beta)&=&\log\int_{-\infty}^{+\infty}\exp\biggl\{-\frac
{U(\beta,\alpha)}{2}\biggr\}\,d\alpha\nonumber\\[-9pt]\\[-9pt]
&&{}-\log\int_{-\infty}^{+\infty}\exp\biggl\{-\frac{U(0,\alpha
)}{2}\biggr\}\,d\alpha\nonumber
\end{eqnarray}
and
%
%e3.24 #&#
\begin{equation}\label{defz}
Z_0(\beta)=\sum_{k=1}^\infty\frac{\beta^2\pi^2k^2X_k^2}{(\pi
^2k^2+\beta^2)\pi^2k^2}+\sum_{k=1}^\infty\log\biggl(\frac{\pi
^2k^2}{\pi^2k^2+\beta^2}\biggr).
\end{equation}
Furthermore, there exists a sequence of local maxima $\hat{\beta
}_{n},\hat{\gamma}_{n}$ of
$L_n(\beta,\gamma)$ converging in distribution to
$\tilde{\beta}_{\mathrm{MLE}},\tilde{\gamma}_{\mathrm{MLE}}$, the global maximum of
the limiting process
$U^*(\beta,\gamma)$. If model (\ref{fmla16}) has, at most, one unit
root, then for the
estimators $\hat c_1$ and~$\hat c_2$, we have
%
%e3.25 #&#
\begin{equation}\label{conv1}
\sqrt{n}\pmatrix{
\hat c_1-c_1 \cr
\hat c_2-c_2}
\dconv\mathrm{N}\biggl(0,\left[\matrix{
1-c_2^2 & c_1(1-c_2) \cr
c_1(1-c_2) & 1-c_2^2}\right]\biggr).
\end{equation}
\end{theorem}
%
%re3.2 #&#
\begin{rem}\label{remp1}
The equivalence in distribution of the processes $U^*(\beta)$ and
$\frac{1}{2}Z_0(\beta)$ is given in Theorem 4.3 in Davis and Song
\cite{lsd1}. As mentioned in Davis and Dunsmiur \cite{dd},
convergence on C(-$\infty$,0] does not necessarily imply convergence
of the corresponding global maximizers. Additional arguments were
required to show that the maximum likelihood estimator converged in
distribution to the global maximizer of the limit process. We
suspect that the same holds here for $\hat{\beta}_{\mathrm{MLE},n}$ and
$\hat{\gamma}_{\mathrm{MLE},n}$ and simulation results, some of which are
contained in Sections \ref{sec4} and \ref{sec5}, bear this out.
\end{rem}
%
%re3.3 #&#
\begin{rem}\label{remp2}
To establish the
convergence in (\ref{conv1}), if there is exactly one unit root, then
\begin{eqnarray*}
\sqrt{n}(\hat{c}_1-c_1)&=&-\frac{\hat{\beta}_{\mathrm{MLE}}}{\sqrt
{n}}-\hat{\gamma}_{\mathrm{MLE}}
\dconv-\tilde{\gamma}_{\mathrm{MLE}}=\frac{N}{\operatorname{var}(N)}\\[-2pt]
&\dequi&
\mathrm{N}(0,1-\alpha_0^2)=\mathrm{N}(0,1-c_2^2),\\[-2pt]
\sqrt{n}(\hat{c}_2-c_2)&=&\hat{\gamma}_{\mathrm{MLE}}+\frac{\alpha_0\hat
{\beta}_{\mathrm{MLE}}}{\sqrt{n}}+\frac{\hat{\gamma}_{\mathrm{MLE}}\hat{\beta
}_{\mathrm{MLE}}}{n}
\dconv\tilde{\gamma}_{\mathrm{MLE}}\\[-2pt]
&=&-\frac{N}{\operatorname{var}(N)}\dequi
\mathrm{N}(0,1-\alpha_0^2)=\mathrm{N}(0,1-c_2^2).
\end{eqnarray*}
Here, we use the fact that $\tilde{\beta}_{\mathrm{MLE}}<\infty$ a.s. as
stated in (Theorem 4.3 in \cite{lsd1}). One can also calculate the
limiting asymptotic covariance of $\hat c_1$ and $\hat c_2$ as
\begin{eqnarray*}
-\operatorname{var}(\tilde{\gamma}_{\mathrm{MLE}})&=&-(1-\alpha_0^2)=
-(1+\alpha_0)(1-\alpha_0)\\
&=&c_1(1-c_2).
\end{eqnarray*}
\end{rem}
%
%re3.4 #&#
\begin{rem}\label{rem1}
The above theorem says that when $|\alpha_0|<1$ and $\theta_0=1$, we
have a similar asymptotic result for $c_1$ and $c_2$ as in
the invertible case. If we only consider the original
parameters $c_1$ and $c_2$, the effect of the unit root
disappears in the limit. But $\sqrt{n}(\hat{c}_1-c_1)$ and
$\sqrt{n}(\hat{c}_2-c_2)$ are perfectly dependent in the
limit, since $c_1(1-c_2)=1-c_2^2$.
\end{rem}
%
%re3.5 #&#
\begin{rem}\label{rem2}
The estimated roots $\hat\theta$ and $\hat\alpha$ calculated from
$\hat c_1$ and
$\hat c_2$ are asymptotically independent.
Interestingly, $\tilde{\beta}_{\mathrm{MLE}}$ corresponding to the unit root
in MA(2) has exactly the same distribution as the
$\tilde{\beta}_{\mathrm{MLE}}$ in the MA(1) case. So the pile-up and other
properties of $\tilde{\beta}_{\mathrm{MLE}}$ follow exactly from those in the
MA(1) case. It may seem surprising that the unit root in
the MA(2) model (when there is only one unit root) behaves
asymptotically just like the unit root in MA(1) case. To see this,
consider the situation where we are given the parameter $\alpha$ and
$\alpha=\alpha_0$. In this case, $\gamma=0$ and
\begin{eqnarray*}
L_n(\beta,0)-L_n(0,0)&\dconv&\log\int_{-\infty}^{+\infty}\exp
\biggl\{-\frac{U(\beta,\eta^*)}{2}\biggr\}\,d\eta^*\\
&&{}-\log\int_{-\infty}^{+\infty}\exp\biggl\{-\frac{U(0,\eta
^*)}{2}\biggr\}\,d\eta^*,
\end{eqnarray*}
which is the limiting process of the exact profile log-likelihood in
the MA(1) case. On the other hand when $\alpha$ is
given, $\theta$ becomes the only parameter that needs to be
estimated
%
%e3.26 #&#
\begin{equation}
X_t=(1-\alpha_0\mathrm{B})(1-\theta\mathrm{B})Z_t.
\end{equation}
Because of the invertibility of the operator $1-\alpha_0\mathrm{B}$,
we can get an intermediate process $Y_t$ by inverting the operator.
Namely,
%
%e3.27 #&#
\begin{equation}
Y_t:=\frac{1}{(1-\alpha_0\mathrm{B})}X_t=\sum_{k=0}^\infty\alpha
_0^kX_{t-k}=(1-\theta\mathrm{B})Z_t.
\end{equation}
Since we are dealing with asymptotics, inverting the operator
$1-\alpha_0\mathrm{B}$ is feasible. Therefore, the transformed
process $Y_t$ is indeed an MA(1) process with the true parameter
$\theta_0=1$. Then it follows naturally that the properties of the
estimator of $\theta$ in this situation should be equivalent to
those of $\theta$ in a unit root MA(1) process.
\end{rem}

%s3.2 #&#
\subsection{MA(2) with two unit roots}
In moving from the unit root problem for the MA(1)\vadjust{\goodbreak} model to the MA(2)
model, several new and challenging
problems arise. In this subsection, we discuss some
issues when there are two unit roots in the MA polynomial.

%s3.2.1 #&#
\subsubsection{\texorpdfstring{Case 2: $c_2=1$ and $c_1\neq\pm2$}{Case 2: c_2=1 and c_1 /= +-2}}

This corresponds to the case that the true parameters are on the
boundary $c_2=1$, that is, the boundary AC in Figure \ref{figtri},
which means the two roots live on the unit circle and are not real
valued. Denote\vspace*{2pt} the two generic complex valued roots of the MA
polynomial by\vspace*{2pt}
$\phi=re^{\vec{i}\theta}$ and $\bar{\phi}=re^{-\vec{i}\theta}$. To
avoid confusion in notation, we use $\vec{i}$ to represent $\sqrt{-1}$.
A rather different representation of the residuals $r_i$ is used
in this case, that is,
%
%e3.28 #&#
\begin{eqnarray}\label{resi1}
r_i&=&\frac{(\phi_0-\phi)(\phi-\bar\phi_0)}{\phi-\bar\phi}\sum
_{j=-1}^{i-1}\phi^{i-j-1}Z_j\nonumber\\[-1pt]
&&{}+\frac{(\bar\phi_0-\bar\phi)(\phi_0-\bar\phi)}{\phi-\bar
\phi}\sum_{j=-1}^{i-1}\bar\phi^{i-j-1}Z_j\nonumber\\[-1pt]
&&{}-\frac{\phi^{i+1}-\bar\phi^{i+1}}{\phi-\bar\phi
}(z_{\mathrm{init},0}-Z_0)\\[-1pt]
&&{}+\frac{\phi\bar\phi(\phi^{i}-\bar\phi
^{i})}{\phi-\bar\phi}(z_{\mathrm{init},-1}-Z_{-1})\nonumber\\[-1pt]
&&{}+\frac{(\phi+\bar\phi-\phi_0-\bar\phi_0)(\phi^{i+1}-\bar\phi
^{i+1})}{\phi-\bar\phi}Z_{-1}.\nonumber
\end{eqnarray}
We also adopt the parameterization for $r$, $\theta$ and two initial
variables given~by
\begin{eqnarray*}
r&=&1+\frac{\beta}{n}\quad\mbox{and}\quad\theta=\theta_0+\frac{\gamma}{n},\\[-1pt]
z_{\mathrm{init},0}&=&Z_0+\frac{\sigma_0\eta_1}{\sqrt{n}}\quad\mbox
{and}\quad z_{\mathrm{init},-1}=Z_{-1}+\frac{\sigma_0\eta_2}{\sqrt{n}}.
\end{eqnarray*}
Again, we study the limiting process of
${-}2\sum_{i=-1}^n\frac{r_iZ_i}{\sigma_0^2}\,{+}\,\sum_{i=-1}^n\frac
{r_i^2}{\sigma_0^2}$.
Here~we only present the first term of
$\sum_{i=-1}^n\!\frac{r_iZ_i}{\sigma_0^2}$ for illustration; the limit
of the~other terms can be derived in a similar fashion. By Theorem
2.8 in \cite{lsd1}, we obtain
\begin{eqnarray*}
&&\frac{1}{n}\sum_{i=0}^n\sum_{j=-1}^{i-1}\phi^{i-j}\frac
{Z_jZ_i}{\sigma_0^2}\\[-1pt]
&&\qquad=\sum_{i=0}^n\Biggl(\sum_{j=-1}^{i-1}\biggl(1+\frac{\beta
}{n}\biggr)^{i-j}\exp\biggl\{\vec{i}\gamma\frac{i-j}{n}\biggr\}
\frac{e^{-\vec{i}\theta_0j}Z_j}{\sqrt{n}\sigma_0}\Biggr)\frac
{e^{\vec{i}\theta_0i}Z_i}{\sqrt{n}\sigma_0}\\[-1pt]
&&\qquad\dconv\int_0^1\int_0^se^{\beta(s-t)+\vec{i}\gamma
(s-t)}\,d\overline{\mathbb{W}}(t)\,d\mathbb{W}(s),
\end{eqnarray*}
where $\mathbb{W}(t)$ is a two-dimensional Brownian motion,
$\mathbb{W}(t)=W_1(t)+\vec{i}W_2(t)$,
$\overline{\mathbb{W}}(t)=W_1(t)-\vec{i}W_2(t)$ and $W_1(t)$ and
$W_2(t)$ are the corresponding weak limits of the sum
\[
W_{1,n}(t)=\sum_{k=0}^{[nt]}\cos(k\theta_0)\frac{Z_k}{\sqrt
{n}\sigma_0} \quad\mbox{and}\quad
W_{2,n}(t)=\sum_{k=0}^{[nt]}\sin(k\theta_0)\frac{Z_k}{\sqrt
{n}\sigma_0}.
\]
The weak convergence of $W_{1,n}(t)$ and $W_{2,n}(t)$ to two
independent Brownian
motions is guaranteed by Theorem 2.2 in Chan and Wei \cite{weichan}.

By Theorem 2.1 in \cite{lsd1} we have
\[
\frac{1}{\sqrt{n}}\sum_{i=0}^n\phi^{i}\frac{Z_i}{\sigma_0}\dconv
\int_0^1e^{\beta
s+\vec{i}\gamma s}\,d\mathbb{W}(s).
\]
Therefore, (\ref{resi1}) leads to
%
%e3.29 #&#
\begin{eqnarray}\label{complexlim}
-2\sum_{i=-1}^n\frac{r_iZ_i}{\sigma_0^2}
&\dconv& 4\Re\biggl\{
(\gamma\cos\theta_0-\gamma\sin\theta_0+\beta\cos\theta_0+\vec
{i}\beta\sin\theta_0)e^{-\vec{i}\theta_0}\nonumber\\
&&\hspace*{54pt}{}\times\int_0^1\int
_0^se^{\beta(s-t)+\vec{i}\gamma(s-t)}\,d\overline{\mathbb
{W}}(t)\,d\mathbb{W}(s)\biggr\}\nonumber\\[-8pt]\\[-8pt]
&&{}+4\eta_1\Re\biggl\{\frac{e^{\vec{i}\theta_0}}{2\vec{i}\sin
\theta_0}\int_0^1e^{\beta
s+\vec{i}\gamma
s}\,d\mathbb{W}(s)\biggr\}\nonumber\\
&&{}-4\eta_2\Re\biggl\{\frac{1}{2\vec{i}\sin\theta_0}\int
_0^1e^{\beta
s+\vec{i}\gamma s}\,d\mathbb{W}(s)\biggr\},\nonumber
\end{eqnarray}
where $\Re\{\cdot\}$ means the real part of a complex function. The
weak limit of $\sum_{i=-1}^nr_i^2/\sigma_0^2$ can also be computed
in an analogous manner using Corollary~2.10 in \cite{lsd1}. However,
the weak limit of $\sum_{i=-1}^nr_i^2/\sigma_0^2$ has an even more
complicated form than (\ref{complexlim}).

By integrating out the auxiliary variables, the exact likelihood can
be recovered as well. However, the form of the joint likelihood
function is much more complicated than the one computed in the one
unit root case. The asymptotic properties and pile-up probabilities
in this case remain unknown.\looseness=1

%s3.2.2 #&#
\subsubsection{Case 3: $c_2=1$ and $c_1=-2$}\label{sec322}

This corresponds to the vertex A in the $\bigtriangledown$-region in
Figure \ref{figtri}. It is convenient to first consider a special
case of local asymptotics when the approach to
the corner is through the boundary $-c_1-c_2=1$.
With this constraint, the dimension of the parameters
has been reduced from two to one. We parameterize the MA(2) in this
case by
%
%e3.30 #&#
\begin{equation}\label{condmodel}
X_t=Z_t-(\theta+1)Z_{t-1}+\theta Z_{t-2}
\end{equation}
and define a $Z_{\mathrm{init}}$ and a $Y_{\mathrm{init}}$ as in (\ref{fmla19}), but
with different normalization, that is,
%
%e3.31 #&#
\begin{equation}\label{ypara}
\theta=1+\frac{\beta}{n},\qquad
Y_{\mathrm{init}}=Y_0+\frac{\sigma_0\eta_1}{n^{{3}/{2}}}
\quad\mbox
{and}\quad
Z_{\mathrm{init}}=Z_{-1}+\frac{\sigma_0\eta_2}{\sqrt{n}}.
\end{equation}
Then, with the help of the theorems in Davis and Song \cite{lsd1},
it follows that
%
%e3.32 #&#
\begin{eqnarray}\label{yinit}
U_n(\beta,\eta_1,\eta_2)&=&-2\sum_{i=-1}^n\frac{r_iZ_i}{\sigma
_0^2}+\sum_{i=-1}^n\frac{r_i^2}{\sigma_0^2}\nonumber\\
&\dconv&2\beta\int_0^1\int_0^se^{\beta(s-t)}\,dW(t)\,dW(s)\nonumber\\[-8pt]\\[-8pt]
&&{}+2\eta_2 \int_0^1e^{\beta
s}\,dW(s)-\frac{2\eta_1}{\beta}\int_0^1(1-e^{\beta
s})\,dW(s)\nonumber\\
&&{}+\int_0^1\biggl(\beta\int_0^se^{\beta(s-t)}\,dW(t)+\eta_2e^{\beta
s}-2\eta_1\frac{1-e^{\beta s}}{\beta}\biggr)^2\,ds.\nonumber
\end{eqnarray}
There is a connection between this limiting process and the one in
(\ref{fmla27}) derived for the limiting process for an MA(1) model with
a nonzero mean. Notice that in~(\ref{fmla27}), $U(\eta_0,\beta,\alpha)$
is exactly the process we just derived with~$\eta_1$ and~$\eta_2$
replaced by $\alpha$ and $\eta_0$. This leads us to an interesting
connection of the mean term in the lower order MA model and the initial
value in the higher order MA model, which we will discuss further in
the Section \ref{sec6}.

Alternatively, if we do not impose the constraint $-c_1-c_2=1$, there
are two possible ways to parameterize the roots. First, the vertex
can be approached through the real region, where
$c_1=-\theta-\alpha$, $c_2=\theta\alpha$ and the roots are
parameterized further as
\[
\theta=1+\frac{\beta}{n} \quad\mbox{and}\quad \alpha=1+\frac
{\gamma}{n},
\]
which makes
\[
c_1=-1-\biggl(1+\frac{\beta+\gamma}{n}\biggr) \quad\mbox{and}\quad
c_2=1+\frac{\beta+\gamma}{n}+o\biggl(\frac{1}{n}\biggr).
\]
The second parameterization is through the complex region, in which the
roots are $re^{\vec{i}\theta}$ and $re^{-\vec{i}\theta}$ with
$c_1=-2r\cos(\theta)$, $c_2=r^2$. The radius and the
angular parts are further parameterized as
\[
r=1+\frac{\beta}{n} \quad\mbox{and}\quad \theta=\frac{\gamma}{n},
\]
which implies
\[
c_1=-1-\biggl(1+\frac{2\beta}{n}\biggr)+o\biggl(\frac{1}{n}
\biggr) \quad\mbox{and}\quad  c_2=1+\frac{2\beta}{n}+o\biggl(\frac
{1}{n}\biggr).
\]
Therefore, in either case, if we ignore the higher order terms,
$c_1$ and $c_2$ can be approximated as
\[
c_1=-1-\biggl(1+\frac{\zeta}{n}\biggr) \quad\mbox{and}\quad  c_2=1+\frac{\zeta}{n}.
\]
This parameterization, however, is exactly the one we have seen
in the conditional case, which suggests that one of the unit roots has
pile-up with probability one asymptotically while the other unit root
behaves like the unit root in the conditional case; see (\ref
{condmodel}) and (\ref{yinit}). This claim is
also supported by the simulation results; see Table \ref{t8} in
Section \ref{sec5}.

%s4 #&#
\section{Testing for a unit root in an MA(2) model}\label{sec4}

A direct application of the results in the previous section is
testing for the presence of a unit root in the MA(2) model. For the
testing problem, we extend the idea of a generalized likelihood
ratio test proposed in Davis, Chen and Dunsmuir \cite{cdd} to the
MA(2) case. Tests based on $\hat{\beta}_{\mathrm{MLE}}$ are also considered
in this section. We will compare these tests with the score-type
test of Tanaka \cite{tpaper}.

To specify our hypothesis testing problem in the MA(2) case, the
null hypothesis is $H_0$: there is exactly one unit root in the MA
polynomial, and the alternative is $H_A$: there are no unit roots.
The asymptotic theory of the previous section allows us to
approximate the nominal power against local alternatives. To set up
the problem, for the model
\[
X_t=Z_t-\biggl(\alpha+1+\frac{\beta}{n}\biggr)Z_{t-1}+\alpha
\biggl(1+\frac{\beta}{n}\biggr)Z_{t-2}
\]
with $|\alpha|<1$. We want to test $H_0\dvtx \beta=0$ versus $H_A\dvtx
\beta<0$.

To describe the test based on the generalized likelihood ratio, let
$\mathrm{GLR}_n=2(L_n(\hat\beta_{\mathrm{MLE}},\hat\gamma_{\mathrm{MLE}})-L_n(0,\hat\gamma
_{\mathrm{MLE},0}))$,
where $\hat\gamma_{\mathrm{MLE},0}$ is the MLE of $\gamma$ when $\beta=0$. An
application\vspace*{2pt} of Theorem \ref{thm4} gives $\mathrm{GLR}_n\dconv
L^*(\tilde\beta_{\mathrm{MLE}},\tilde\gamma_{\mathrm{MLE}})-L^*(0,\tilde\gamma
_{\mathrm{MLE}})=U^*(\tilde\beta_{\mathrm{MLE}})$,
where $L^*(\beta,\gamma)$ and $U^*(\beta)$ are\vspace*{1pt} given
in~(\ref{newadded})\break
and~(\ref{defu}) and~$\tilde\gamma_{\mathrm{MLE}}=-N/\operatorname{var}(N)$. Notice that the limit
distribution of $\mathrm{GLR}_n$ only depends on $\tilde\beta_{\mathrm{MLE}}$, and
$\gamma$ serves as a nuisance parameter, which does not play a role
in the limit. Define the $(1-\alpha)$th asymptotic quantile~$b_{\mathrm{GLR}}(\alpha)$ and
 $b_{\mathrm{MLE}}(\alpha)$ as
\[
\mathbf{P}\bigl(U^*(\tilde\beta_{\mathrm{MLE}})>b_{\mathrm{GLR}}(\alpha)
\bigr)=\alpha \quad\mbox{and}\quad \mathbf{P}\bigl(\tilde{\beta
}_{\mathrm{MLE}}>b_{\mathrm{MLE}}(\alpha)\bigr)=\alpha.
\]

Since the limiting random variables $U^*(\tilde\beta_{\mathrm{MLE}})$ and
$\tilde{\beta}_{\mathrm{MLE}}$ are the same as in the MA(1) unit root case,
the critical values of $b_{\mathrm{GLR}}(\alpha)$ and $b_{\mathrm{MLE}}(\alpha)$ are
the same as those provided in Table 3.2 of Davis, Chen and Dunsmuir
\cite{cdd}.

There has been limited research on the testing for a unit root in
the MA(2) case. One approach, proposed by Tanaka, was based
on a score type of statistic, which is locally best invariant and
unbiased (LBIU). However, implementation of this test requires
choosing a sequence $l_n\rightarrow\infty$ at a~suitable rate. One
choice is $l_n=o(n^{1/4})$, yet this may not always work well,
especially if $\alpha>0$; see also \cite{tpaper}. Next we compare
the power curves of the three tests for sample size $n=50$.

Figure \ref{figf1} below shows the power curves based on MLE, GLR
and LBIU tests, when the invertible root $\alpha$ in the MA(2) model
is $-0.3$ and $-0.5$, respectively. Since the score-type test of Tanaka
is demonstrated to be locally best invariant unbiased, it has a very
small edge on the GLR test up to the local alternative 4 or so.
Thereafter, the GLR test increasingly outperforms the LBIU test by a
wide margin. When the sample size is 50, the local alternative
parameter corresponds to $\theta=1-4/50=0.92$. Also, as seen in Figure~\ref{figf1},
the power function based on the MLE dominates the
power function of the LBIU test for local alternatives greater than 8
or 9.

%f2 #&#
\begin{figure}

\includegraphics{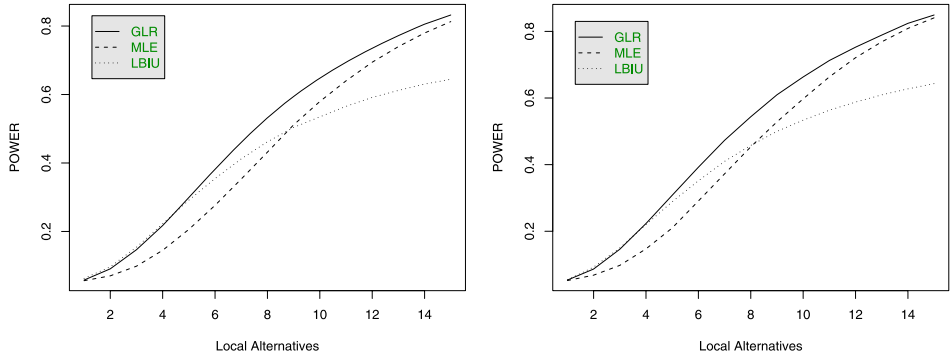}

\caption{Power curve with respect to local alternatives when
$\alpha=-0.3$ (upper) and when $\alpha=-0.5$ (lower). Sample size
$n=50$. The size of the test is set to be 0.05.}\label{figf1}
\end{figure}

In the case when $\alpha>0$ especially for small sample sizes like
50, the behavior of the tests based on MLE and LBIU are very poor.
This is because when $\alpha>0$ and there is one unit root, the two
parameters $c_1$ and $c_2$ lie on the boundary
$-c_1-c_2=1$ which is close to the complex region boundary
$c_1^2-4c_2=0$. But our asymptotic results are derived in
a way which assumes that the two roots are only approaching the
limit through the real region. This holds asymptotically, but in
finite sample cases, when we maximize the likelihood jointly over
$c_1$ and $c_2$, it is likely that the two maximizers
would fall into the complex region. As $\alpha$ gets closer to $-1$
this effect becomes more severe. Thus we do not recommend using the
test based on the MLE when the invertible root is likely to be
negative. Using the test based on MLE usually gives larger size of
the test. The LBIU is not good in this case either as pointed out in
Tanaka \cite{tpaper}. The upper tail probabilities are greatly
underestimated when $\alpha$ gets closer to $-1$, and hence $H_0$ tends to
be accepted much more often. Simulation results show that when the
sample size is 50, and the true $\alpha$ is 0.3 and 0.5, the
corresponding size of the LBIU test is 0.0119 and 0.0015 which are
much smaller than the nominal size 0.05. GLR seems to be the best
among the three choices. This is due to the fact that the GLR only
considers the maximum value of the likelihood ratio instead of the
MLE of $c_1$ and $c_2$. Therefore, even if $\hat c_1$
and $\hat c_2$ are in the complex region, the GLR test can still
be carried out whereas the test based on $\hat\beta_{\mathrm{MLE}}$ is not
even well defined in this case. Although the size of the GLR test is
often slightly greater than the nominal size, GLR gives the best
performance under this situation.

%f3 #&#
\begin{figure}

\includegraphics{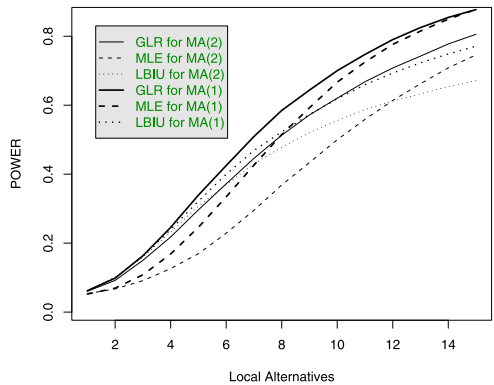}

\caption{Power curve with respect to local alternatives when
$\alpha=0$. Sample size $n=50$. The size of the test is set to be
0.05.}\label{figf2}
\end{figure}

Finally, we compare these tests when $\alpha=0$; that is, the model is in
fact a~unit root MA(1). The test developed for
the MA(2) case is still applicable. The results are summarized in
Figure \ref{figf2}. Clearly, the power functions of the tests
designed for the MA(1) dominate the power functions of their
counterparts designed for the MA(2). However, it is surprising that
for large local alternatives (greater than 9 or so), the GLR for the
MA(2) model outperforms the LBIU for the MA(1) model.

%s5 #&#
\section{Numerical simulations}\label{sec5}

In this section, we present simulation results that illustrate
the theory from Section \ref{sec3}. Realizations were simulated from the
MA(2) process given by
%
%e5.1 #&#
\begin{equation}
X_t=Z_t-(1+\alpha)Z_{t-1}+\alpha Z_{t-2},
\end{equation}
where $\alpha$ takes the values 0.3, 0 and $-0.3$, respectively. The
MA(2) model was replicated 10,000 times for each choice of $\alpha$,
and then
the MLEs for the MA(2) coefficients $\theta_1$ and $\theta_2$ were
calculated for each replicate. The empirical pile-up
probability, the empirical variance and MSE of the MLEs are reported in
Tables \ref{t3} to \ref{t5}.
%
%t1 #&#
\begin{table}
\caption{Summary of the case: $\alpha=0.3$}
\label{t3}
\begin{tabular*}{\tablewidth}{@{\extracolsep{\fill}}lcccccc@{}}
\hline
\textbf{Sample} & \textbf{Pile-up} & \textbf{Variance} & \textbf{MSE}
& \textbf{Variance} & \textbf{MSE} & \textbf{Correlation} \\
\textbf{size} & \textbf{probability} & \textbf{of} $\bolds{c_1}$
& \textbf{of} $\bolds{c_1}$ & \textbf{of} $\bolds{c_2}$
& \textbf{of} $\bolds{c_2}$ & \textbf{of} $\bolds{c_1}$
\textbf{and} $\bolds{c_2}$ \\
\hline
\hphantom{00,}25 & 0.5436 & 2.1701 & 2.1970 & 2.4455 & 2.6536 & 0.9347\\ % inserting
%body of the table
\hphantom{00,}50 & 0.6041 & 1.4063 & 1.4118 & 1.4967 & 1.5553 & 0.9644\\
\hphantom{0,}100 & 0.6234 & 1.1108 & 1.1108 & 1.1490 & 1.1636 & 0.9815\\
\hphantom{0,}400 & 0.6398 & 0.9788 & 0.9788 & 0.9854 & 0.9890 & 0.9953\\
1,000 & 0.6437 & 0.9290 & 0.9290 & 0.9327 & 0.9338 & 0.9981\\  %
%[1ex] adds vertical space
\hline%inserts single line
\end{tabular*}
%
% is used to refer this table in the text
\end{table}
Notice that the numbers in the tables for the variance and the MSE
are reported for the normalized estimates $\sqrt{n}(\hat c_i-c_i)$, $i=1,2$.

%t2 #&#
\begin{table}
\caption{Summary of the case: $\alpha=0$ [MA(1) with a unit root]} %
\label{t4}
\begin{tabular*}{\tablewidth}{@{\extracolsep{\fill}}lcccccc@{}}
\hline
\textbf{Sample} & \textbf{Pile-up} & \textbf{Variance} & \textbf{MSE}
& \textbf{Variance} & \textbf{MSE} & \textbf{Correlation} \\
\textbf{size} & \textbf{probability} & \textbf{of} $\bolds{c_1}$
& \textbf{of} $\bolds{c_1}$ & \textbf{of} $\bolds{c_2}$
& \textbf{of} $\bolds{c_2}$ & \textbf{of}
$\bolds{c_1}$ \textbf{and} $\bolds{c_2}$ \\
\hline
\hphantom{00,}25 & 0.5870 & 2.1624 & 2.1629 & 2.5037 & 2.6355 & 0.8792\\
\hphantom{00,}50 & 0.6182 & 1.3661 & 1.3670 & 1.4690 & 1.5053 & 0.9378\\
\hphantom{0,}100 & 0.6220 & 1.1661 & 1.1670 & 1.2082 & 1.2224 & 0.9662\\
\hphantom{0,}400 & 0.6318 & 1.0440 & 1.0441 & 1.0544 & 1.0578 & 0.9918\\
1,000 & 0.6334 & 1.0329 & 1.0330 & 1.0351 & 1.0384 & 0.9966\\
\hline
\end{tabular*}
\end{table}
%
%t3 #&#
\begin{table}
\caption{Summary of the case: $\alpha=-0.3$}
\label{t5}
\begin{tabular*}{\tablewidth}{@{\extracolsep{\fill}}lcccccc@{}}
\hline
\textbf{Sample} & \textbf{Pile-up} & \textbf{Variance} & \textbf{MSE}
& \textbf{Variance} & \textbf{MSE} & \textbf{Correlation} \\
\textbf{size} & \textbf{probability} & \textbf{of} $\bolds{c_1}$
& \textbf{of} $\bolds{c_1}$ & \textbf{of} $\bolds{c_2}$
& \textbf{of} $\bolds{c_2}$ & \textbf{of}
$\bolds{c_1}$ \textbf{and} $\bolds{c_2}$ \\
\hline
\hphantom{00,}25 & 0.6171 & 1.8370 & 1.8806 & 2.1654 & 2.2287 & 0.7950\\
\hphantom{00,}50 & 0.6347 & 1.2820 & 1.3053 & 1.3647 & 1.3820 & 0.8938\\
\hphantom{0,}100 & 0.6447 & 1.0748 & 1.0853 & 1.1215 & 1.1299 & 0.9397\\
\hphantom{0,}400 & 0.6472 & 0.9245 & 0.9267 & 0.9316 & 0.9339 & 0.9822\\
1,000 & 0.6511 & 0.9232 & 0.9242 & 0.9256 & 0.9263 & 0.9933\\  %
\hline
\end{tabular*}
\end{table}

As seen in the tables, the correlation of $\hat c_1$ and
$\hat c_2$ is increasing to 1 with the sample size. The
variances and the MSEs are converging to the theoretical value
$1-c_2^2$. As pointed out in \cite{dd} and \cite{cdd}, the asymptotic results
work remarkably well even for small sample sizes in the MA(1) case. Here,
although the pile-up probability is still 0.6518, the rates vary
depending on~$\alpha$. For $\alpha>0$, rates are slow while for
$\alpha<0$ rates are much faster. From the derivation of the asymptotic
results, there are error terms in the likelihood that vanish
asymptotically and contribute to a more lethargic
rate of convergence. Again the
asymptotic results were derived assuming the roots are always in the
real region, which only holds asymptotically. When the sample size
is small and $\alpha>0$, the MLEs of $c_1$ and $c_2$ are
more likely to be in the complex region than those when $\alpha<0$.
Thus the limiting process would approximate the likelihood function
poorly when $\alpha>0$, which in turn results in less pile-up in
smaller sample sizes.

Table \ref{t8} summarizes the pile-up effects for the model
considered in Section~\ref{sec322}, where the two roots of the MA
polynomial are both 1. In one realization, the estimators are said
to exhibit a pile-up if the MLEs of $c_1$ and $c_2$ are on the
boundary $-c_1-c_2=1$.

As seen in the table, the pile-up probability is increasing to 1
with sample size. However, the claimed 100\% probability of pile-up
is not a good approximation for small sample sizes. Even when
$n=500$, the pile-up is only about 80\%.

%
%t4 #&#
\begin{table}
\tablewidth=205pt
\caption{Pile-up probabilities for the case: $c_1=-2$}
\label{t8}
\begin{tabular*}{\tablewidth}{@{\extracolsep{\fill}}lc@{}}
\hline
\textbf{Sample size} & \textbf{Pile-up probability}\\
\hline
\hphantom{0,}100 & 0.246 \\
\hphantom{0,}500 & 0.804 \\
1,000 & 0.961 \\
5,000 & 0.999 \\
\hline
\end{tabular*}
\end{table}
%

%s6 #&#
\section{Unit roots and differencing}\label{sec6}

As pointed out in Section \ref{sec322}, there is a link between the mean
term in the lower order MA model and the initial value in the higher
order MA model. To illustrate this, consider the simple case when
\[
Y_t=\mu_0+Z_t,
\]
where $\{Z_t\}\sim$ i.i.d. $(0,\sigma_0^2)$. So $Y_t$ is an i.i.d.
sequence with a common mean. It is
clear that
\[
\sqrt{n}(\hat{\mu}-\mu_0)\dconv\mathrm{N}(0,\sigma_0^2),
\]
where $\hat{\mu}$ is the MLE of $\mu$ obtained by maximizing the
objective Gaussian likelihood function. Now suppose we difference
the time series to obtain
\[
X_t=(1-\mathrm{B})Y_t=Z_t-Z_{t-1},
\]
which becomes an MA(1) process with a unit root. The initial
value as defined before of this differenced process is
\[
Z_{\mathrm{init}}=Z_0=Y_0-\mu_0.
\]
From the results in Theorem 4.2 in \cite{lsd1}, if it is known that an
MA(1) time series has a unit root, that is, $\beta=0$, we have
\[
U(\beta=0,\alpha)=2\alpha W(1)+\alpha^2.
\]
Clearly, $\tilde{\alpha}=-W(1)$ and with our parameterization of
$z_{\mathrm{init}}$, we have
\begin{eqnarray*}
\hat{\alpha}&=&\frac{\sqrt{n}(z_{\mathrm{init}}-Z_0)}{\sigma_0}=\frac
{\sqrt{n}(Y_0-\hat{\mu}-Y_0+\mu_0)}{\sigma_0}\\
&=&-\frac{\sqrt{n}(\hat{\mu}-\mu_0)}{\sigma_0}\dconv\tilde
\alpha=-W(1)\dequi\mathrm{N}(0,1),
\end{eqnarray*}
which is consistent with the classical result. Therefore we can
conclude that whenever we have an MA model with a unit root, the
information stored in the initial value comes from the information
of the mean term from the undifferenced series. So differencing the
series will not get rid of the mean parameter; instead, differencing
creates a new parameter $Z_{\mathrm{init}}$ which behaves like the mean in
the undifferenced series and its effect persists even
asymptotically. With this, we can now explain easily the result in
(\ref{fmla28}). Turning to a little more
complicated model consisting of i.i.d. noise and a linear trend, that is,
%
%e6.1 #&#
\begin{equation}\label{fmla29}
Y_t=\mu_0+b_0t+Z_t,
\end{equation}
which, after differencing, delivers an MA(1) model with a~unit root and
a~nonzero mean given by
\[
X_t=(1-\mathrm{B})Y_t=b_0+Z_t-Z_{t-1}.
\]
From (\ref{fmla28}), we know
$n^{{3}/{2}}(\hat{b}-b_0)\dconv\mathrm{N}(0,12\sigma_0^2)$. But
this can be obtained much more easily by analyzing the model
(\ref{fmla29}). This is just a simple application of linear
regression, and we can get exactly the same asymptotic result for
$\hat b$.

Now consider the model from Section \ref{sec2},
\[
Y_t=b_0+Z_t-\theta Z_{t-1},
\]
where $\theta=1+\frac{\beta}{n}$ is near or on the unit circle. By
differencing we obtain
\[
X_t=(1-\mathrm{B})Y_t=Z_t-(1+\theta)Z_{t-1}+\theta Z_{t-2}.
\]
If we define $Z_{\mathrm{init}}$ as before and
\[
Y_{\mathrm{init}}=Y_0=b_0+Z_0-Z_{-1},
\]
then $y_{\mathrm{init}}-Y_{\mathrm{init}}$ can be viewed as $\hat b-b_0$. Since $\hat b$
converges at the rate of~$n^{3/2}$, so does $y_{\mathrm{init}}$. This
explains the parametrization given in (\ref{ypara}) as well as the
resemblance of (\ref{fmla27}) and (\ref{yinit}).

%s7 #&#
\section{Going beyond second order}\label{sec7}

The techniques proposed in this paper can be adapted to handle the
unit root problem for MA($q$) with $q\geq3$. However, the complexity
of the argument, mostly in terms of bookkeeping, also increases with
the order $q$. In this section, we outline the procedure for the
MA(3) case, from which extensions to larger orders are
straightforward.\looseness=1\vadjust{\goodbreak}

Suppose $\{X_t\}$ follows an MA(3) model, which is parameterized in
terms of the reciprocals of the zeros of the MA polynomial, that is,
%
%e7.1 #&#
\begin{eqnarray}\label{ma3model}
X_t&=&Z_t-(\theta_0+\phi_0+\psi_0)Z_{t-1}\nonumber\\
&&{} +(\theta_0\phi_0+\theta_0\psi_0+\phi_0\psi
_0)Z_{t-2}-\theta_0\phi_0\psi_0Z_{t-3} \nonumber\\
&=&(1-\theta_0B)(1-\phi_0B)(1-\psi_0B)Z_t\\
&=&(1-\theta_0B)(1-\phi_0B)Y_t\nonumber\\
&=&(1-\theta_0B)W_t.\nonumber
\end{eqnarray}
For simplicity, assume $\theta_0\neq\phi_0\neq\psi_0$. Now we form
two intermediate processes~$Y_t$ and~$W_t$ and consider three
augmented initial variables defined by $Z_{\mathrm{init}}=Z_{-2}$,
$Y_{\mathrm{init}}=Z_{-1}+\psi_0Z_{\mathrm{init}}$ and $W_{\mathrm{init}}=Y_0+\phi_0Y_{\mathrm{init}}$.
Similar arguments as in Section \ref{sec3} show that the joint likelihood of
$(\mathbf{X},W_{\mathrm{init}},Y_{\mathrm{init}},Z_{\mathrm{init}})$ has a~simple form given by
\[
f_{\mathbf{X},W_{\mathrm{init}},Y_{\mathrm{init}},Z_{\mathrm{init}}}(\mathbf
{x}_n,w_{\mathrm{init}},y_{\mathrm{init}},z_{\mathrm{init}})=\prod_{j=-2}^nf_Z(z_j).
\]
As in the MA(1) and MA(2) cases, maximizing this joint likelihood is
essentially equivalent to minimizing the objective function
\[
U_n=\frac{1}{\sigma_0^2}\sum_{i=-2}^n(z_i^2-Z_i^2).
\]
The key to this analysis is to write out the explicit expression for
$z_i$ which is basically an estimator for $Z_i$. The following
equations are straightforward to derive:
%
%e7.2 #&#
%e7.3 #&#
%e7.4 #&#
\begin{eqnarray}
\label{wk}
w_k&=&\sum_{l=1}^k\theta^{k-l}X_l+\theta^kw_{\mathrm{init}},\\
\label{yj}
y_j&=&\sum_{k=1}^j\phi^{j-k}w_k+\phi^jw_{\mathrm{init}}+\phi
^{j+1}y_{\mathrm{init}},\\
\label{zi}
z_i&=&\sum_{j=1}^i\psi^{i-j}y_j+\psi^iw_{\mathrm{init}}+\psi^i(\phi+\psi
)y_{\mathrm{init}}+\psi^{i+2}z_{\mathrm{init}}.
\end{eqnarray}
Plugging (\ref{wk}) into (\ref{yj}), we obtain
%
%e7.5 #&#
\begin{equation}\label{yjnew}
y_j=\sum_{k=1}^j\frac{\theta^{j-k+1}-\phi^{j-k+1}}{\theta-\phi
}X_k+\frac{\theta^{j+1}-\phi^{j+1}}{\theta-\phi}w_{\mathrm{init}}+\phi
^{j+1}y_{\mathrm{init}},
\end{equation}
and plugging this into (\ref{zi}), we obtain
\begin{eqnarray*}
z_i&=&\sum_{j=1}^i\fracb{\theta\phi(\theta^{i-j+1}-\phi
^{i-j+1})\\
&&\hspace*{17.2pt}{}+\theta\psi(\psi^{i-j+1}-\theta^{i-j+1})+\phi\psi(\phi
^{i-j+1}-\psi^{i-j+1})}\\
&&\hspace*{13pt}{}\times\bigl({(\theta-\psi)(\psi-\phi)(\phi-\theta
)}\bigr)^{-1}X_j\\
&&{}+\biggl(\frac{\theta^2(\theta^i-\psi^i)}{(\theta-\phi)(\theta
-\psi)}-\frac{\phi^2(\phi^i-\psi^i)}{(\theta-\phi)(\phi-\psi
)}+\psi^i\biggr)w_{\mathrm{init}}\\
&&{}+\frac{\phi^{i+2}-\psi^{i+2}}{\phi-\psi
}y_{\mathrm{init}}+\psi^{i+2}z_{\mathrm{init}}.
\end{eqnarray*}
While this is a more complicated looking expression than the one
encountered in the MA(2) case, the coefficient of $X_j$ in the sum
looks very similar to~(\ref{fmla4th}), only with more terms. Now
replacing $X_j$ with (\ref{ma3model}), $z_i$ can be written as
%
%e7.6 #&#
\begin{eqnarray}\label{generalform}
z_i&=&Z_i-\sum_{j=-2}^{i-1}C^z_{i,j}Z_j\nonumber\\
&&{}-C^w_i(w_{\mathrm{init}}-W_{\mathrm{init}})-C^y_i(y_{\mathrm{init}}-Y_{\mathrm{init}})-C^z_i(z_{\mathrm{init}}-Z_{\mathrm{init}})
\\
&=&Z_i-r_i,\nonumber
\end{eqnarray}
where $C^z_{i,j}$ is the coefficient for $Z_j$ in $z_i$ and is a
combination of $\theta^{i-j}$, $\phi^{i-j}$ and~$\psi^{i-j}$, and
$C^w_i$, $C^y_i$ and $C^z_i$ are coefficients for
$w_{\mathrm{init}}-W_{\mathrm{init}}$, $y_{\mathrm{init}}-Y_{\mathrm{init}}$ and $z_{\mathrm{init}}-Z_{\mathrm{init}}$.
They are linear combinations of $\theta^{i}$, $\phi^{i}$ and
$\psi^{i}$. For illustration, assume the MA(3) model has only one
unit root with $|\psi_0|<1$, $|\phi_0|<1$ and $\theta_0=1$. We can
then reparameterize the parameters as
\[
\theta=1+\frac{\beta}{n},\qquad \beta\leq0,\qquad \phi=\phi
_0+\frac{\alpha}{\sqrt{n}}\quad\mbox{and}\quad\psi=\psi_0+\frac{\gamma
}{\sqrt{n}},
\]
and the initial values as
\[
w_{\mathrm{init}}=W_{\mathrm{init}}+\frac{\sigma_0\eta_w}{\sqrt{n}},\qquad
y_{\mathrm{init}}=Y_{\mathrm{init}}+\frac{\sigma_0\eta_y}{\sqrt{n}}
\quad\mbox{and}\quad
z_{\mathrm{init}}=Z_{\mathrm{init}}+\frac{\sigma_0\eta_z}{\sqrt{n}}.
\]
Then the objective function $U_n$ becomes
%
%e7.7 #&#
\begin{eqnarray}\label{genlik}
&&U_n(\beta,\alpha,\gamma,\eta_w,\eta_y,\eta_z)\nonumber\\
&&\qquad=-2\sum
_{i=-2}^n\frac{r_iZ_i}{\sigma_0^2}+\sum_{i=-2}^n\frac{r_i^2}{\sigma
_0^2}\nonumber\\[-8pt]\\[-8pt]
&&\qquad=-2\sum_{i=-2}^n\Biggl(\sum_{j=-2}^{i-1}C_{i,j}^z\frac{Z_j}{\sigma
_0}+\frac{C_i^w\eta_w}{\sqrt{n}}+\frac{C_i^y\eta_y}{\sqrt
{n}}+\frac{C_i^z\eta_z}{\sqrt{n}}\Biggr)\frac{Z_i}{\sigma
_0}\nonumber\\
&&\qquad\quad{}+\sum_{i=-2}^n\Biggl(\sum_{j=-2}^{i-1}C_{i,j}^z\frac{Z_j}{\sigma
_0}+\frac{C_i^w\eta_w}{\sqrt{n}}+\frac{C_i^y\eta_y}{\sqrt
{n}}+\frac{C_i^z\eta_z}{\sqrt{n}}\Biggr)^2.\nonumber
\end{eqnarray}
Because of the special structure of $C^z_{i,j}$, $C_i^w$, $C_i^y$
and $C_i^z$, the sum in (\ref{genlik}) consists of terms that have a
similar structure to quantities like
\[
\frac{1}{n}\sum_{i=-2}^n\sum_{j=-2}^{i-1}\biggl(1+\frac{\beta
}{n}\biggr)^{i-j}\frac{Z_j}{\sigma_0}\frac{Z_i}{\sigma_0}
\quad\mbox{and}\quad\sum_{i=-2}^n\biggl(1+\frac{\beta}{n}\biggr)^{i}\frac
{Z_i}{\sqrt{n}\sigma_0},
\]
that were used in the MA(1) and MA(2) cases. By using a martingale
central limit theorem and theorems proved in Davis and Song
\cite{lsd1}, one can establish the weak convergence\vspace*{1pt}
of~$U_n(\beta,\alpha,\gamma,\eta_w,\eta_y,\eta_z)$ to a random
element~$U(\beta,\alpha,\gamma,\eta_w,\eta_y,\eta_z)$ in
$\mathbb{C}(\mathbb{R}^6)$. Now arguing as in Section \ref{sec3}, the initial
variables can be integrated out, and the limiting process of the
exact profile log-likelihood can be established.

For general $q>3$, the residual $r_i=z_i-Z_i$ has the form
\[
r_i=\sum_{j=-q+1}^{i-1}C_{i,j}^zZ_j+\sum_{k=1}^qC_i^{k}(\mathrm{init}_k-\mathrm{INIT}_k),
\]
where $\{\mathrm{INIT}_1,\ldots,\mathrm{INIT}_q\}$ are $q$ augmented initial
variables, defined either through the i.i.d. random variables $Z_t$ or
through the intermediate processes like $Y_t$ in the above example.
Furthermore,\vspace*{1pt} $C_{i,j}^z$ is only a linear combination of
$(\theta_1^{i-j},\ldots,\theta_q^{i-j})$, where
$(\theta_1,\ldots,\theta_q)$ are reciprocals of the roots of the
MA(q) polynomial. Coefficients $C_{i}^k$, $k=1,\ldots,q$, are only
linear combinations of
$(\theta_1^{i},\ldots,\theta_q^{i})$. This special
structure of $r_i$ allows us to apply the weak convergence theorems
proved in Davis and Song \cite{lsd1} to find the limiting process of
$U_n=-2\sum_{i=-q+1}^n r_iZ_i/\sigma_0^2+\sum_{i=-q+1}^n
r_i^2/\sigma_0^2$, from which the limiting behavior of the maximum
likelihood estimators of the $\theta_i$'s can be derived.

\section*{Acknowledgments}

We would like to thank the referees and the Associate
Editor for their insightful comments, which were incorporated into
the final version of this paper.

%suskaldyti doi

% imsref loaded by lrinkeviciute, 2011-11-29 11:18:14

\printaddresses

\end{document}